\documentclass[a4paper,12pt,leqno]{amsart}
\usepackage{amsmath,amssymb,mathrsfs}
\usepackage[matrix,arrow]{xy} 
\usepackage{xcolor,hyperref}
\hypersetup{colorlinks=true,citecolor=black,linkcolor=black}
\newcommand{\h}{\hbox}
\newcommand{\q}{\quad}
\newcommand{\nin}{\noindent}
\newcommand{\bs}{\par\bigskip}
\newcommand{\ms}{\par\medskip}
\newcommand{\sk}{\par\smallskip}

\newcommand{\msn}{\par\medskip\noindent}
\newcommand{\skn}{\par\smallskip\noindent}
\newcommand{\ges}{\geqslant}
\newcommand{\les}{\leqslant}
\newcommand{\one}{\hskip1pt}
\newcommand{\mopl}{\hbox{$\bigoplus$}}
\newcommand{\mcup}{\hbox{$\bigcup$}}
\newcommand{\msum}{\hbox{$\sum$}}
\newcommand{\mprod}{\hbox{$\prod$}}
\newcommand{\A}{{\mathscr A}}

\newcommand{\D}{{\mathscr D}}
\newcommand{\F}{{\mathscr F}}
\newcommand{\G}{{\mathscr G}}
\newcommand{\Hc}{{\mathscr H}}

\newcommand{\Lc}{{\mathscr L}}
\newcommand{\M}{{\mathscr M}}
\newcommand{\OO}{{\mathscr O}}
\newcommand{\Tc}{{\mathscr T}}
\newcommand{\Fb}{{\mathscr F}^{\ssb}}
\newcommand{\Fsb}{{\mathscr F}_{\!\!s}^{\ssb}}
\newcommand{\Lt}{\widetilde{L}}

\newcommand{\Yt}{\widetilde{Y}}
\newcommand{\Zt}{\widetilde{Z}}
\newcommand{\rht}{\widetilde{\rho}}
\newcommand{\sit}{\widetilde{\si}}
\newcommand{\DD}{{\mathbb D}}

\newcommand{\PP}{{\mathbb P}}
\newcommand{\R}{{\mathbb R}}

\newcommand{\Q}{{\mathbb Q}}
\newcommand{\C}{{\mathbb C}}
\newcommand{\N}{{\mathbb N}}
\newcommand{\Z}{{\mathbb Z}}
\newcommand{\mm}{{\mathfrak m}}
\newcommand{\rt}{\widetilde{r}}

\newcommand{\al}{\alpha}
\newcommand{\ga}{\gamma}
\newcommand{\la}{\lambda}
\newcommand{\si}{\sigma}
\newcommand{\ep}{\varepsilon}

\newcommand{\ddd}{{\rm d}}
\newcommand{\Gr}{{\rm Gr}}
\newcommand{\Om}{\Omega}
\newcommand{\Si}{\Sigma}
\newcommand{\Sing}{{\rm Sing}}

\newcommand{\dX}{\delta\!\one _X}
\newcommand{\DR}{{\rm DR}}
\newcommand{\MHM}{{\rm MHM}}

\newcommand{\bl}{\bigl}
\newcommand{\br}{\bigr}
\newcommand{\ssb}{\raise.15ex\h{${\scriptscriptstyle\bullet}$}}
\newcommand{\ssc}{\one\raise.15ex\hbox{${\scriptstyle\circ}$}\one}
\newcommand{\onto}{\twoheadrightarrow}
\newcommand{\into}{\hookrightarrow}
\newcommand{\simto}{\,\,\rlap{\hskip1.3mm\raise1.4mm\hbox{$\sim$}}\hbox{$\longrightarrow$}\,\,}
\newcommand{\pl}{\one{+}\one}
\newcommand{\mi}{\one{-}\one}
\newcommand{\ins}{\,{\in}\,}
\newcommand{\tos}{\,{\to}\,}
\newcommand{\eq}{\,{=}\,}
\newcommand{\nes}{\,{\ne}\,}
\newcommand{\defs}{\,{:=}\,}
\newcommand{\gess}{\,{\ges}\,}
\newcommand{\less}{\,{\les}\,}
\newcommand{\slt}{\,{<}\,}
\newcommand{\sgt}{\,{>}\,}

\newcommand{\sst}{\,{\subset}\,}
\newcommand{\stm}{\,{\setminus}\,}
\newcommand{\col}{\,{:}\,}
\newcommand{\pv}{\par\noindent\hskip5mm\verb}
\renewcommand{\qedsymbol}{}
\makeatletter
\renewcommand\section{\@startsection{section}{1}{0pt}{-3.5ex plus -1ex minus -.2ex}{2.3ex plus.2ex}{\centering\normalfont\bfseries}}
\makeatother
\theoremstyle{plain}
\newtheorem{thm}{Theorem}[section]
\newtheorem{cor}[thm]{Corollary}
\newtheorem{prop}[thm]{Proposition}
\newtheorem{lem}[thm]{Lemma}
\newtheorem{ithm}{Theorem}
\newtheorem{icor}{Corollary}
\newtheorem{iprop}{Proposition}

\theoremstyle{definition}
\newtheorem{rem}[thm]{Remark}
\newtheorem{conj}[thm]{Conjecture}
\newtheorem{exam}[thm]{Example}

\newtheorem*{claim}{Claim}

\newtheorem{iques}{Question}

\begin{document}
\title[$L$-classes of projective varieties]{Constant coefficient and intersection complex $L$-classes of projective varieties}
\author{J. Fern\'andez de Bobadilla}
\address{Javier Fern\'andez de Bobadilla:
(1) IKERBASQUE, Basque Foundation for Science, Maria Diaz de Haro 3, 48013,
Bilbao, Basque Country, Spain;
(2) BCAM, Basque Center for Applied Mathematics, Mazarredo 14, 48009 Bilbao,
Basque Country, Spain;
(3) Academic Colaborator at UPV/EHU.}
\email{jbobadilla@bcamath.org}
\author{I. Pallar\'es}
\address{Irma Pallar\'es Torres: University of Cantabria, Department of Mathematics, Statistics and Computation, Avda. de los Castros s/n. 39005 Santander, Spain}
\email{irma.pallares@unican.es}
\author{M. Saito}
\address{Morihiko Saito: RIMS Kyoto University, Kyoto 606-8502 Japan}
\email{msaito@kurims.kyoto-u.ac.jp}
\thanks{J.F.B. was supported by the Basque Government through the BERC 2022-2025 program, by the Ministry of Science and Innovation: BCAM Severo Ochoa accreditation CEX2021-001142-S / MICIN / AEI / 10.13039/501100011033. and by the Spanish Ministry of Science, Innovation and Universities, project reference PID2020-114750GB-C33.
I.P. was supported by 12B1423N from FWO, Research Foundation Flanders.
M.S. was partially supported by JSPS Kakenhi 15K04816.}
\begin{abstract} For a projective variety $X$, we have the intersection complex $L$-classes $L_*(X)$ defined by Goresky-MacPerson using cohomotopy and also the constant coefficient $L$-class $L^c_*(X)$ defined by applying an $L$-class transformation (or $T_{1*}$) to a cubic hyperresolution of $X$. These coincide if $X$ is a $\mathbb Q$-homology manifold. We show that the two $L$-classes $L_*(X)$ and $L^c_*(X)$ differ if they do by replacing $X$ with an intersection of general hyperplane sections which has only $\mathbb Q$-homologically isolated singularities. Finding a good sufficient condition for the non-coincidence of $L_*(X)$ and $L^c_*(X)$ is thus reduced to the latter case, where a necessary and sufficient condition has been obtained in terms of the Hodge signatures of stalks of intersection complex in our previous paper. In the case of projective hypersurfaces having only isolated singularities, the difference between $L_*(X)$ and $L^c_*(X)$ is given by the Hodge signatures of the link cohomologies at singular points, and the Hodge signatures of the vanishing cohomologies give the difference between $L^c_*(X)$ and the virtual $L$-class of $X$, that is, the image by a retraction map of the $L$-class of a smooth deformation of $X$ in an ambient smooth projective variety $Y$ in the very ample case.
\end{abstract}
\maketitle

\section*{Introduction}
The Hirzebruch characteristic (cohomology) class of a compact complex manifold $X$, denoted by $T_y^*(X)\ins H^{2\ssb}(X,\Q)[y]$, was introduced in \cite{Hi}. It specializes to the Chern, Todd, and $L$-classes of $X$ at $y\eq{-}1,0$ and $1$ respectively. This is generalized to the theory of Hirzebruch characteristic (homology) classes for singular varieties $X$, which are denoted by $T_{y*}(X)$, and belong to $H^{\rm BM}_{2\ssb}(X,\Q)[y]$, the Borel-Moore homology of $X$ with even degrees tensored with $\Q[y]$, see \cite{BSY}.
\sk
Specializing $T_{y*}$ at $y\eq{-}1$, the class transformation $T_{-1*}$ gives the scalar extension of MacPherson's Chern class transformation $c_*{\otimes}\Q$ in \cite{Ma}. At $y\eq 0$, $T_{0*}(X)$ coincides with the Todd class $td_*(X)$ in \cite{BFM} (see also \cite{Fu}) if $X$ has only Du Bois singularities (so that $\underline{\Omega}^0_X=\OO_X$). At $y\eq 1$, assuming $X$ is {\it compact,} we have the commutative diagram
\begin{equation}\begin{gathered}
\xymatrix{& K_0\bl(\MHM(X,\R)\br) \ar[d]^{\rm Pol} \\
K_0({\rm Var}/X) \ar[ru]^{\rm Hdg} \ar[r]^{sd} \ar[rd]_{T_{1*}}
& \Om_{\R}(X) \ar[d]^{L_*}\\
& H_{2\one\ssb}(X,\Q) }\end{gathered}
\label{1} \end{equation}
whose lower and upper parts are constructed in \cite{BSY} (reducing to the smooth projective case \cite{Hi}, see also Sections\,\,\ref{S2} and \ref{S3} below) and in \cite{FPS} respectively. Here Hdg, Pol, and $sd\one $ stand for Hodge, polarization, and self-duality respectively, $\Om_{\R}(X)$ is the cobordism class group of self-dual bounded $\R$-complexes with constructible cohomology sheaves on $X$, and $L_*$ is a homology $L$-class transformation extending \cite{GM}, see \cite{CS}, \cite{Ba1}, \cite{Yo}, \cite{BSY}, and also Section\,\,\ref{S3} below. Here the $L$-class transformation $L_*$ is given by $L_k\defs(-1)^{k(k-1)/2}\Lt_k$ with $\Lt_*$ a sign modified $L$-class transformation defined in Section~\ref{S3}. Note that self-pairings are used in the definition of $\Om_{\R}(X)$ in this paper (see also \cite{FPS}), and this is compatible with the usual definition as in \cite{BSY} according to \cite{AMS} (where duality isomorphisms are mainly used rather than self-pairings and a question noted at the end of \cite[Remark 1.1e]{FPS} may be avoided).
\sk
As remarked in \cite{BSY}, the image of $[X]\defs[X\,{\buildrel{\!\!\rm id}\over\to}\,X]\ins K_0({\rm Var}/X)$ by $sd$ in \eqref{1} does not necessarily coincide with the cobordism class of the intersection complex $[({\rm IC}_X\R,S_X)]$ (even up to sign) unless $X$ is smooth, where $S_X$ is the natural pairing. So the specialization at $y\eq1$ of the Hirzebruch class $T_{y*}(X)\,({=}\,T_{y*}([X])\eq T_{y*}([\R_{h,X}]))$ may be different in general from the {\it topological intersection complex $\one L$-class\one} $L_*(X)$, that is, the homology $L$-class in \cite{GM}, which is defined by applying $L_*$ to $[({\rm IC}_X\R,S_X)]$. Here ${\rm Hdg}([X])\eq[\R_{h,X}]$, and $\R_{h,X}$ denotes a bounded complex of mixed $\R$-Hodge modules such that its underlying $\R$-complex is the constant sheaf $\R_X$ and the zeroth cohomology $H^0(X,\Q_{h,X})$ is a pure Hodge structure of type $(0,0)$.
Note that the specialization $T_{1*}(X)$ can be defined by using a cubic hyperresolution of $X$ (see \cite[3.2]{mhc}, \cite[Section 2]{FP}), and coincides by the commutativity of \eqref{1} with the {\it constant coefficient\one} $L$-class $L^c_*(X)$ defined by applying $L_*$ to the image of $[X]$ by $sd\one $.
\sk
It has been conjectured that the two $L$-classes $L^c_*(X)$ and $L_*(X)$ coincide in the case $X$ is a $\Q$-{\it homology manifold\one} (that is, if ${\rm IC}_X\R\eq\R_X[d_X]$ with $d_X\defs\dim X$), see \cite{BSY}. This conjecture has been shown in certain special cases by many people (see for instance \cite{Ba2}, \cite{CMSS1}, \cite{CMSS2}, \cite{MS}), and is recently proved in \cite{FP} and \cite{FPS} for the projective and general compact cases respectively; more precisely it is shown there that the image of $[X\,{\buildrel{\rm id\,}\over{\to}}\,X]$ in $\Om_{\R}(X)$ by the transformation $sd_{\R}$ coincides with $[({\rm IC}_X\R,S_X)]$ for projective and general compact $\Q$-homology manifolds respectively.
\sk
It is easy to see that the {\it converse\one} of the above conjecture about the coincidence of $L^c_*(X)$ and $L_*(X)$ does not necessarily hold. We are then interested in the following.

\begin{iques} \label{Q1}
Is there a good sufficient condition for the non-coincidence of $L^c_*(X)$ and $L_*(X)\,$?
\end{iques}

A good necessary condition for non-coincidence may be that $X$ is not a $\Q$-homology manifold. In order to make a comparison, we may sometimes need to assume that the class of the intersection complex Hodge module ${\rm IC}_X\R_h$ belongs to the image of the morphism Hdg in \eqref{1}. This condition is satisfied if $X$ is not very complicated, for instance in the isolated singularity case (see Remark\,\,\ref{R1.1} below), but the general case seems quite non-trivial, since the stability of the image of Hdg by taking {\it direct factors\one} (or Ker, Coker of morphisms) of pure Hodge modules seems apparently false, see Section\,\,\ref{S10} below and also \cite{Sc2}.
\sk
We say that a hypersurface $X'\sst X$ is {\it stratified-smooth\one} if its scheme-theoretic intersection with any stratum of a Whitney stratification $\mathscr S$ of $X$ is a (reduced) smooth hypersurface of the stratum (in the case $X'$ is locally the restriction of a hypersurface $Y'$ in a local ambient smooth variety $Y$ containing locally $X$, the intersection of $Y'$ with each stratum of $\mathscr S$ is transversal, and $X'$ does not meet any 0-dimensional strata), see Remark\,\,\ref{R1.2} below.
\sk
We say that the $L$-classes $L_*(X)$ and $L^c_*(X)$ {\it differ on $X$ at degree\one} $2\dX$ if their components in $H_{2\dX}(X,\Q)$ differ. Here $\dX\defs\dim\Si_X$ with $\Si_X$ the {\it non-$\Q$-homology-manifold locus\one} of $X$, that is, the complement of the largest open subvariety of $X$ which is a $\Q$-homology manifold. This may be called the $\Q$-{\it homological singular locus.}
\sk
In the projective case we say that the two $L$-classes $L_*(X)$ and $L^c_*(X)$ {\it strongly differ at degree\one} $2\dX$ (with projective embedding fixed) if the degree $2\dX$ part of the images of the two $L$-classes in the homology group of an ambient projective space differ. 
\sk
In this paper we prove the following.

\begin{iprop} \label{P1}
For a compact variety $X$, the constant coefficient $L$-class $L^c_*(X)$ and the intersection complex $L$-class $L_*(X)$ differ at degree $2\dX$ if they do after replacing $X$ with a stratified-smooth hypersurface $X'\sst X$ $($for instance, a general hyperplane section in the projective case$)$. A similar assertion holds with ``differ" replaced by ``strongly differ" in the case $X$ is projective.
\end{iprop}

Since the normal bundle of $X'$ in $X$ is nontrivial, we need the topological filtration $G$ (see Section~\ref{S4} below) in an essential way. Note that there is a counterexample to the converse of the first assertion of Proposition\,\,\ref{P1} using Theorem\,\,\ref{T1} just below, see Remark\,\,\ref{R6.1}.
\begin{ithm} \label{T1}
For a projective variety $X$, the degree $2\dX$ part of the difference between the images of the constant coefficient and intersection complex $L$-classes in the homology of an ambient projective space is given by the sum of the reduced modified Hodge signature $\sit_x$ in Theorem\,\,{\rm\ref{T2}} below, which is applied to each connected component of the intersection of sufficiently general $\dX$ hyperplane sections of $X$, where a projective embedding of $X$ is fixed.
\end{ithm}

\begin{icor} \label{C1}
For a projective variety $X$ with $\dX\gess1$, the two $L$-classes $L_*(X)$ and $L^c_*(X)$ strongly differ at degree $2\dX$ if and only if it holds for a general hyperplane section of $X$.
\end{icor}

Here $X$ is {\it not\one} assumed to be pure-dimensional. The intersection complexes of irreducible components of $X$ are shifted depending on their dimensions as in the formula \eqref{2} below. The proofs of Proposition\,\,\ref{P1} and Theorem\,\,\ref{T1} are not trivial even with the definition of intersection complex $L$-class. Indeed, the Todd class transformation is not compatible with pullbacks, and we need the {\it topological filtration\one} on the Grothendieck group of coherent sheaves in an essential way. We consider the ``strong difference" because of the counterexample to the converse of the first assertion of Proposition\,\,\ref{P1}, see Remark~\ref{R6.1} below. This also gives a counterexample to Corollary\,\,\ref{C1} with ``strongly" omitted.
\sk
Question\,\,\ref{Q1} for a projective variety is thus reduced to the isolated singularity case. Here we have a necessary and sufficient condition in terms of the mixed Hodge structures on the stalks of the intersection complexes of the irreducible components $X_i$ of $X$ (which are obtained by using the theory of mixed Hodge modules \cite{mhp}, \cite{mhm}, \cite{mhc}) as follows:

\begin{ithm}[{\cite[Theorem 2]{FPS}}] \label{T2}
Assume $X$ is a connected compact variety with $\dX\eq0$. Then the two $L$-classes $L_*(X)$ and $L^c_*(X)$ coincide if and only if
\begin{equation}
\aligned&\msum_{x\in\Si_X}\,\sit_x=0\q\q\h{with}\\&\sit_x:=\msum_{j,p,k\in\Z}\,(-1)^{j+p}\dim_{\C}\Gr_F^p\Gr^W_{2k}H^j_x-1,\\&H^j_x\defs\mopl_{X_i\ni x}\,\Hc^j({\rm IC}_{X_i}\C[-d_{X_i}])_x,\endaligned
\label{2} \end{equation}
where the $X_i$ are irreducible components of $X$.
\end{ithm}

Concerning the definition of $\sit_x$ in \eqref{2}, note that there are many {\it cancellations of signs,} since $i\one ^{q-p}\eq i\one ^{p-q}\eq(-1)^{k-p}$ if $p\pl q\eq 2k$, and $(-1)^{w(w+1)/2}\eq(-1)^k$ if $w\eq 2k$, see also a remark before \cite[Theorem 1]{FPS}. We can {\it omit\one} $\Gr^W_{2k}$ and the summation over $k$ in the definition of the {\it reduced modified Hodge signature\one} $\sit_x$ so that the formula agrees with the specialization at $y\eq1$ of the definition of Hirzebruch classes in \eqref{1.1} below (with $X\eq pt$). Indeed, the contribution of the {\it odd weight\one} part vanish by the {\it Hodge symmetry.} Note also that if $x\,{\notin}\,\Si_X$, we have $\dim H^j_x\eq\delta_{j,0}$ and $\sit_x\eq0$. In the proof of \cite[Theorem 2]{FPS}, we assumed that $\dim\Sing\,X\eq0$, but the same argument applies to the case $\dX\eq0$ using Proposition\,\,\ref{P3.4}.
\sk
We now assume $X$ is a {\it very ample reduced divisor\one} on a smooth projective variety $Y$ having only (geometrically) {\it isolated singularities,} for instance, $Y\eq\PP^{d_X+1}$. Let $Z$ be a sufficiently general member of the complete linear system defined by $X$ (where $Z$ is smooth, since $X$ is very ample). There is a one-parameter family $g\col\Yt\tos\PP^1$ such that $g^{-1}(0)\eq X$, $g^{-1}(\infty)=Z$, and $\Yt$ is the blowing-up of $Y$ along $X\cap Z$. We denote by $\varphi_{g,1}\R_{h,\Yt}[d_X]$ the unipotent monodromy part of the vanishing cycle Hodge module $\varphi_g\R_{h,\Yt}[d_X]$, which is identified with the direct sum of vanishing cohomologies $V_x\defs H^{d_X}(F_{\!g,x},\R)$ with $F_{\!g,x}$ the Milnor fiber around $x\ins{\rm Sing}\,X$.
Let $L_{X,x}\defs S_x\cap X$, the link of $X$ at $x\ins{\rm Sing}\,X$, where $S_x$ is a sufficiently small sphere around $x$ in $Y$.
\sk
Set $N\defs\log T_u$ with $T\eq T_sT_u$ the Jordan decomposition of the monodromy. We have the decomposition $V_x\eq V_{x,1}\oplus V_{x,\ne1}$ in a compatible way with the action of the monodromy $T$ so that the eigenvalues of the action of $T$ on $V_{x,1}$ and $V_{x,\ne1}$ are respectively 1 and not 1. The weight filtration $W$ on $V_{x,1}$, $V_{x,\ne1}$ is given by the monodromy filtration shifted by $d_Y\eq d_X{+}1$ and $d_X$ respectively.
Let $P_N\Gr^W_{d_X+1+k}V_{x,1}$, $P_N\Gr^W_{d_X+k}V_{x,\ne1}$ be the $N$-primitive part of the $W$-graded quotients of $V_{x,1}$, $V_{x,\ne1}$ for $k\gess0$. They vanish otherwise, see also Section\,\,\ref{S8} below. For $j\ins\Z$, define
\begin{equation*}
\si^j_{x,1}\defs\si(P_N\Gr^W_jV_{x,1}),\q\si^j_{x,\ne1}\defs\si(P_N\Gr^W_jV_{x,\ne1}),
\end{equation*}
where $\si(H)\defs\msum_{p\in\Z}\,(-1)^p\dim_{\C}\Gr_F^pH_{\C}$ denotes the Hodge signature of a Hodge structure $H$, which is the specialization at $p\eq1$ of $T_{y*}(H)\eq\chi_y(H)\defs\msum_{p\in\Z}\,(-y)^p\dim_{\C}\Gr_F^pH_{\C}$. This coincides with the signature of ${\rm Pol}(H)$, since the sign $(-1)^{w(w+1)/2}$ appears in the definition of Pol (see \cite[Theorem 3]{FPS}) where $w$ is the weight of Hodge structure, see \cite[p.~427]{Sc} and also Lemma~\ref{L3.6} below. We have the following.

\begin{ithm} \label{T3}
Let $X$ be a very ample reduced divisor on a connected smooth projective variety $Y$ having only $($geometrically$)$ isolated singularities. Let $L_{X,x},V_x,Z$ be as above. For any $x\ins{\rm Sing}\,X$ and $k\gess0$, there are isomorphisms of mixed Hodge structures
\begin{equation}
\Gr^W_{d_X-1-k}H^{d_X-1}(L_{X,x},\R)=(P_N\Gr^W_{d_X+1+k}V_{x,1})(k{+}1),
\label{3} \end{equation}
together with the equalities
\begin{equation}
L_j(X)=L^c_j(X)=r_*L_j(Z)\,\,\,(\forall\,j>0),
\label{4} \end{equation}
\vskip-5mm
\begin{equation}
\aligned L^c_0(X)\mi L_0(X)&=\msum_{x\in{\rm Sing}\,X,\,k\ges0}\,(-1)^{k+1}\si_{x,1}^{d_X+1+k},\\ L_0(Z)\mi L^c_0(X)&=\msum_{x\in{\rm Sing}\,X,\,k\ges0}\,\bl(\si_{x,1}^{d_X+1+2k}\pl\si_{x,\ne1}^{d_X+2k}\br),\\ L_0(Z)\mi L_0(X)&=\msum_{x\in{\rm Sing}\,X,\,k\ges0}\,\bl(\si_{x,1}^{d_X+2+2k}\pl\si_{x,\ne1}^{d_X+2k}\br).\endaligned
\label{5} \end{equation}
\end{ithm}

Here $L_*(Z)$ is the homology $L$-class of $Z$, and $H_0(X,\Z)$, $H_0(Z,\Z)$ are identified with $\Z$. The direct image $r_*L_{\ssb}(Z)$ by the retraction $r$ is independent of the choice of a path $\ga$ on the parameter space of hypersurfaces of $Y$, see Section\,\,\ref{S8} below for the retraction $r$.
Note that the non-unipotent vanishing cohomologies $V_{x,\ne1}$ do not necessarily vanish even if $X$ is a $\Q$-homology manifold as in Brieskorn's examples of exotic spheres, and the $\si^j_{x,1}$, $\si^j_{x,\ne1}$ {\it vanish\one} if $j$ is {\it odd.} It does not seem easy to generalize the argument to the non-isolated singularity case using \cite{MSY} (for instance there does not seem to exist a morphism from a smooth hypersurface of the line bundle whose general fibers are hypersurfaces of $Y$).
\sk
Finally we note the following.

\begin{iprop} \label{P2}
There is a projective variety with two $L$-classes $L_*(X)$ and $L^c_*(X)$ strongly differing, although the constant coefficient and intersection cohomologies coincide $($where ${\rm Sing}\,X\cong\PP^1)$.
\end{iprop}

This implies that the difference between the two $L$-classes cannot be determined by the constant coefficient and intersection cohomologies. (So ``intersection cohomology $L$-class" is not a good terminology.) Note that it is easy to get examples where the two $L$-classes {\it coincide,} although the constant coefficient and intersection cohomologies {\it differ.} Indeed, the {\it odd\one} weight part does {\it not\one} contribute to the difference between the two $L$-classes in \eqref{2}.
\sk
We thank J.\,Sch\"urmann and the referee for valuable comments about this paper.

\tableofcontents
\numberwithin{equation}{section}

\section{Hirzebruch classes of Hodge modules} \label{S1}
Let $\MHM(X,\R)$ be the category of mixed $\R$-Hodge modules on a complex variety $X$. We use $\R$-Hodge modules since the transformation Pol cannot be defined for the Grothendieck group of $\Q$-Hodge modules, see \cite[Section 2.5]{FPS}.
For $\M\ins D^b\MHM(X)$, its {\it homology Hirzebruch characteristic class\one} is defined by
\begin{equation}
\aligned T_{y*}(\M)&:=td_{(1+y)*}\bl(\DR_y[\M]\br)\in H^{\rm BM}_{2\ssb}(X,\Q)\bl[y,\h{$\frac{1}{y(y+1)}$}\br]\q\q\h{with}\\ \DR_y[\M]&:=\msum_{j,p\in\Z}\,(-1)^j\,\bl[\Hc^j\Gr_F^p\DR(\M)\br]\,(-y)^p\in K_0(X)[y,y^{-1}].\raise11pt\h{ }\endaligned
\label{1.1} \end{equation}
Here $F^p=F_{-p}$, and we define
\begin{equation}
td_{(1+y)*}:K_0(X)[y,y^{-1}]\to H^{\rm BM}_{2\ssb}(X,\Q)\bl[y,\h{$\frac{1}{y(y+1)}$}\br]
\label{1.2} \end{equation}
to be the {\it scalar extension\one} of the Todd class transformation
\begin{equation*}
td_*:K_0(X)\to H^{\rm BM}_{2\ssb}(X,\Q).
\end{equation*}
(denoted by $\tau$ in \cite{BFM}) followed by {\it multiplication by\one} $(1{+}y)^{-k}$ on the degree $2k$ part, see \cite{BSY}. Actually, by \cite[Prop.~5.21]{Sc} we have
\begin{equation}
T_{y*}(\M)\in H^{\rm BM}_{2\ssb}(X,\Q)[y,y^{-1}].
\label{1.3} \end{equation}
The de Rham complex $\DR(\M)$ is defined by using local embeddings of $X$ into smooth varieties. It is shifted by the dimension of smooth varieties as usual. The graded quotient cohomology sheaves $\Hc^i\Gr_F^p\DR(\M)$ are then independent of local embeddings, and are $\OO_X$-modules by \cite[Lemma 3.2.6]{mhp}. Note that the Todd class transformation $td_*$ is compatible with the pushforward by proper morphisms.
\sk
The {\it homology Hirzebruch characteristic class} $T_{y*}(X)$ of a complex algebraic variety $X$ is defined by applying the above definition to the case $\M\eq\R_{h,X}\defs a_X^*\R_h$, where $\R_h$ is the trivial Hodge structure of type $(0,0)$ with $a_X\col X\tos pt$ the structure morphism (see \cite{mhm} for $a_X^*$), that is,
\begin{equation}
\aligned &T_{y*}(X):=T_{y*}(\R_{h,X})=td_{(1+y)*}\DR_y[X]\in H^{\rm BM}_{2\ssb}(X,\Q)[y],\\&\h{with}\q\q\q\DR_y[X]:=\DR_y[\Q_{h,X}].\endaligned
\label{1.4} \end{equation}
This coincides with the definition using the Du Bois complex by \cite{mhc}. It is shown that $T_{y*}(X)$ belongs to $H^{\rm BM}_{2\ssb}(X,\Q)[y]$ in \cite{BSY}.

\begin{rem} \label{R1.1}
The morphism Hdg in \eqref{1} is defined by assigning $[f_!\R_{h,X'}]\ins K_0(\MHM(X,\R))$ to $[X'\,{\buildrel{f}\over\to}\,X]\ins K_0({\rm Var}/X)$. Recall that $\R$-Hodge structures are uniquely determined by their Hodge numbers, and for an elliptic curve $E$, the $\R$-Hodge structure $H^1(E,\R)$ (whose class in the Grothendieck ring coincides with the image of $[\PP^1]\mi[E]$ by Hdg) is independent of $E$. Using the K\"unneth formula, we see that any $\R$-Hodge structure with level at most its weight belongs to the image of Hdg in \eqref{1}. (This is used implicitly in \cite{FPS}.) Here the {\it level\one} of a Hodge structure $(H,F)$ is the difference between the maximal and minimal numbers $p$ with $\Gr_F^pH_{\C}\nes0$. Note also that a negative Tate twist $(-k)$ of Hodge structures corresponds to the product with $[\PP^k]\mi[\PP^{k-1}]$ in $K_0({\rm Var}/pt)$ for $k\sgt0$.
\end{rem}

\begin{rem} \label{R1.2}
The condition for {\it stratified-smooth hypersurface\one} in the introduction implies that the restriction to a stratum $S$ of the differential of a local defining function of $X'$ in an ambient space does not vanish at any intersection point $x$, since the {\it embedding dimension,} that is, $\dim{\mm}_{S,x}/{\mm}_{S,x}^2$, must decrease, where ${\mm}_{S,x}\sst\OO_{S,x}$ is the maximal ideal. Recall that the {\it scheme-theoretic intersection\one} is defined by the sum of the ideals of subvarieties.
\end{rem}

\begin{rem} \label{R1.3}
It is quite well known to specialists that $[X]\ins K_0({\rm Var}/X)$ is represented by $\msum_{l\in\N}\,(-1)^l[X_l]$ with $\ep\col X_{\ssb}\tos X$ a cubic hyperresolution \cite{GNPP} (where the morphisms {\it between\one} the $K_l$ can be neglected in the $K$-group). This can be verified by induction on $\dim X$. If we take a desingularization $\pi\col X'\tos X$ and put $X''\defs\pi^{-1}(X_{\rm sing})$, then we have $[X]\eq[X']\pl[X_{\rm sing}]\mi[X'']$, and the assertion follows if $X''$ and $X_{\rm sing}$ are smooth (for instance if $X$ is a projective cone of a smooth projective variety). In general we apply the inductive hypothesis to cubic hyperresolutions of $X''$ and $X_{\rm sing}$ compatible with $\pi|_{X''}$. Note that the ``mapping cone" of a cubic hyperresolution $X_{\ssb}\tos X$ in general is the single complex of some $n$-ple complex; consider for instance the case $X''$, $X_{\rm sing}$ are smooth, where $n\eq2$, see also \cite[Section 2]{FP}. Here it does not seem necessary to use the assertions that the weight filtration is induced by the truncations $X_{\ssb\ges -k}$ for $-k\ins\N$ in the $X$ compact case and this gives the $E_1$-complex of the weight spectral sequence, see \cite{GNPP} and also \cite[2.3]{mhc}.
\end{rem}

\section{Cohomology L-classes} \label{S2}
For an oriented compact differentiable (that is, $C^{\infty})$ manifold $M$ with $\dim M\eq2n$, the {\it Pontryagin classes\one} $p_i(TM)$ ($i\ins\N$) of the tangent bundle $TM$ are defined by
\begin{equation}
p_i(TM):=(-1)^ic_{2i}(TM{\otimes}\C)\ins H^{4i}(M,\Q).
\label{2.1} \end{equation}
Here the $c_i$ denote the Chern classes, and $TM{\otimes}\C$ is the complexification of the real vector bundle $TM$. The cohomology $L$-{\it class\one} $L^{\ssb}(M)$ of $M$ is defined by
\begin{equation}
L^{(k)}\bl(p_1(TM),\dots,p_{[n/2]}(TM)\br)\ins H^{2k}(M,\Q)\,\,\,\,(k\ins\N),
\label{2.2} \end{equation}
where the $L^{(k)}$ denote the {\it multiplicative sequence of polynomials\one} associated with the analytic function $x/\!\tanh(x)$. The latter is a power series of $x^2$, and $L^{(k)}\eq0$ for $k$ odd, see \cite{Hi}. (The relation to the homology $L$-class $L_*(M)$ is given by the equality $L_k(M)\eq L^{2n-k}(M)\cap[M]$ using Poincar\'e duality.) We have the well known {\it Hirzebruch signature theorem} as follows:

\begin{thm}[see {\cite[Theorem 8.2.2]{Hi}, \cite{Hi2}}] \label{T2.1}
For an oriented compact differentiable manifold $M$ of dimension $4m$, its signature $\si(M)$ is expressed as
\begin{equation}
\si(M)=\bl\langle L^{(2m)}\bl(p_1(TM),\dots,p_m(TM)\br),[M]\br\rangle.
\label{2.3} \end{equation}
\end{thm}

\begin{rem} \label{R2.1}
The construction of the $L$-class in \cite{GM}, \cite{CS} for the shifted constant sheaf $\R_X[d_X]$ in the $X$ smooth compact case is equivalent to the above one using the Hirzebruch signature theorem. Indeed, we may assume that a {\it cohomotopy\one} $\rho\col M\tos S^{2k}$ is given by a differentiable map employing a smooth approximation theorem as in \cite{Ba1}. Since the {\it normal bundle\one} of a general fiber $F$ of $\rho$ is {\it trivial,} we get the equality
\begin{equation}
\aligned\si(F)&=\bl\langle L^{(2m)}\bl(p_1(TF),\dots,p_m(TF)\br),[F]\br\rangle\\&=\bl\langle L^{(2m)}\bl(p_1(TM),\dots,p_m(TM)\br)|_{F},[F]\br\rangle,\endaligned
\label{2.4} \end{equation}
where $d_X{-}k$ is even and $m\defs(d_X{-}k)/2$. This gives the reason why the two constructions of $L$-classes are equivalent (using Poincar\'e duality) at least in the case the homology group $H_{4m}(M,\Q)$ is generated by general fibers of differentiable cohomotopies, see \cite{GM}, \cite{CS}, \cite{Ba1} for a more precise argument. (In a certain case we have to replace $X$ by its product with a sphere of sufficiently high dimension.)
\end{rem}

\begin{rem} \label{R2.2}
For a compact complex manifold $X$ of dimension $n$, the {\it Pontryagin classes\one} $p_i(TM)$ of the real tangent bundle $TM$ of the underlying real differentiable manifold $M$ of $X$ is related to the {\it Chern classes\one} $c_j(TX)$ of the complex tangent bundle $TX$ as follows:
\begin{equation}
\msum_{i\ges0}\,(-1)^ip_i(TM)=\bl(\msum_{j\ges0}\,c_j(TX)\br)\bl(\msum_{j\ges0}\,(-1)^jc_j(TX)\br),
\label{2.5} \end{equation}
see \cite[Theorem 4.6.1]{Hi}. (Note that the real vector bundle $TM$ is identified with the underlying real vector bundle of the complex vector bundle $TX$, and the complexification of $TM$ is isomorphic to the direct sum of $TX$ and its complex conjugation.)
\sk
This implies that the ``formal Pontryagin roots" $y_k$ of $TM$ are given by the $x_l^2$ with $x_l$ the formal Chern roots of the complex vector bundle $TX$ (using the {\it splitting principle\one} if necessary), see for instance \cite[a remark after (15)]{Hi2}. Indeed, by \eqref{2.5} we have
\begin{equation}
\mprod_{k=1}^n\,(1\mi y_k)=\mprod_{l=1}^n\,(1\pl x_l)\cdot\mprod_{l=1}^n\,(1\mi x_l)=\mprod_{l=1}^n\,(1\mi x_l^2).
\label{2.6} \end{equation}
\sk
By the definition of multiplicative sequences \cite{Hi}, this implies that the cohomology $L$-class of $X$ is given by
\begin{equation}
\aligned L^{2\ssb}(X)=\mprod_{l=1}^n\,H(x_l^2)\,\,\,\,\,\h{in}\,\,\, H^{4\ssb}(X,\Q)\q\q\q\q\q\h{with }\\H(z):=\sqrt{z}/\!\tanh(\sqrt{z})=\bl(\msum_{k\ges0}\,z^k/(2k)!\br)\big/\bl(\msum_{k\ges0}\,z^k/(2k{+}1)!\br),\endaligned
\label{2.7} \end{equation}
where $H(z)$ is called the {\it Hirzebruch function.} This argument is generalized to the case of holomorphic vector bundles on compact complex manifolds. Here the power series expansion of the Hirzebruch function is as follows (using Macaulay2 or Singular):
\begin{equation}
\aligned&\h{$H(z)\eq1\pl\tfrac{1}{3}z\mi\tfrac{1}{45}z^2\pl\tfrac{2}{945}z^3\mi\tfrac{1}{4725}z^4\pl\tfrac{2}{93555}z^5\mi\tfrac{1382}{638512875}z^6\pl\tfrac{4}{18243225}z^7\mi\tfrac{3617}{162820783125}z^8$}\\&\h{$\pl\tfrac{87734}{38979295480125}z^9\mi\tfrac{349222}{1531329465290625}z^{10}\pl\tfrac{310732}{13447856940643125}z^{11}\mi\tfrac{472728182}{201919571963756521875}z^{12}\pl\cdots$}\endaligned
\label{2.8} \end{equation}
This is compatible with a table of Bernoulli numbers in \cite[p.\,12]{Hi}.
\end{rem}

\begin{exam} \label{E2.1}
In the case $X\eq\PP^n$, we have the Euler sequence
\begin{equation}
0\tos\OO_X\tos\buildrel{n+1}\over\mopl\OO_X(1)\tos\Tc_X\tos0,
\label{2.9} \end{equation}
where $\Tc_X$ denotes the tangent sheaf. Using the generalization of the above argument on the formal roots of Pontryagin classes and Chern classes (see a remark after \eqref{2.7}), this implies the following.
\end{exam}
\begin{lem} \label{L2.1}
The cohomology $L$-class $L^{2\ssb}(\PP^n)\in H^{4\ssb}(\PP^n,\Q)$ is identified with
\begin{equation}
H(z)^{n+1}\,\,\,{\rm mod}\,\,(z^{[n/2]+1})\in\Q[z]/(z^{[n/2]+1})\eq H^{4\ssb}(\PP^n,\Q),
\label{2.10} \end{equation}
where $z\eq x^2$ with $x\ins H^2(\PP^n,\Q)$ the hyperplane class.
\end{lem}
Applying the {\it Hirzebruch signature theorem\one} (that is, Theorem\,\,\ref{T2.1}) to the case $n\eq2m$ with $m\ins\Z_{>0}$, we get the following

\begin{cor} \label{C2.1}
For any positive integer $m$, the coefficient of $z^m$ in $H(z)^{2m+1}$ is $1$.
\end{cor}
For instance, if $m\eq12$, the cohomology $L$-class of $\PP^{24}$ is given by
\begin{equation}
\aligned&\h{$H(z)^{25}\,\,{\rm mod}\,\,(z^{13})=1\pl\tfrac{25}{3}z\pl\tfrac{295}{9}z^2\pl\tfrac{5090}{63}z^3\pl\tfrac{79211}{567}z^4\pl\tfrac{480296}{2673}z^5\pl\tfrac{1513618772}{8513505}z^6$}\\&\h{$\pl\tfrac{704727224}{5108103}z^7\pl\tfrac{236526822118}{2791213425}z^8\pl\tfrac{72144676811}{1749912255}z^9\pl\tfrac{46326137919619}{2940824761875}z^{10}\pl\tfrac{1073136102266}{231905038365}z^{11}\pl z^{12}.$}\endaligned
\label{2.11} \end{equation}
Here the coefficients of $z^j$ in $H(z)^{n+1}$ are identified with $L^{2j}(\PP^n)\ins H^{4j}(\PP^n,\Q)$ for any $j\less n/2$, but these are not necessarily integers, since we have to divide the ``signature" by the ``degree" of a general fiber of a ``cohomotopy". One may propose the following.

\begin{conj} \label{Con2.1}
The positivity $L^{2j}(\PP^n)\sgt0$ holds for any $j\ins[0,n/2]$.
\end{conj}

\begin{exam} \label{E2.2}
The above argument can be extended to the {\it smooth projective complete intersection\one} case as below (see \cite{MSS0} for the Hirzebruch characteristic classes of singular complete intersections).
\end{exam}

\begin{prop} \label{P2.1}
For a smooth complete intersection $X$ of multidegree $(d_1,\dots,d_r)$ in $\PP^{2m+r}$, the cohomology $L$-class $L^{2\ssb}(X)$ belongs to the non-primitive cohomology $H^{4\ssb}_{\rm nprim}(X,\Q)$ which is the image of $H^{\ssb}(\PP^{2m+r},\Q)$ by definition and is identified with $\Q[z]/(z^{m+1})$, and the $L$-class is expressed as
\begin{equation}
H(z)^{2m+r+1}\one\mprod_{i=1}^r\,H(d_i^{\one 2}z)^{-1}\,\,\,\,\,{\rm mod}\,\,\,(z^{m+1}),
\label{2.12} \end{equation}
with $x,z$ as in {\rm\eqref{2.10}}. Hence the signature of $X$ is given by the coefficient of $z^m$ in this power series multiplied by the degree of $X$ $($which is equal to $\mprod_{i=1}^r\,d_i)$. In particular the signature is the value of a polynomial of $d_1^2,\dots,d_r^2$ of degree at most $2m$ multiplied by $d_1\cdots d_r$.
\end{prop}

\begin{proof}
This follows from the same argument as above using the exact sequence of vector bundles on $X$\,:
\begin{equation*}
0\tos TX\tos TY|_X\tos N_{X/Y}\tos 0,
\end{equation*}
with $Y\defs\PP^{2m+r}$. Indeed, the normal bundle $N_{X/Y}$ of $X$ in $Y$ is isomorphic to the direct sum of the restrictions of the normal bundles $N_{D_i/Y}$ over $X\sst D_i$, where $D_i\defs\{f_i\eq0\}\sst Y$ with $f_i\ins H^0(Y,\OO_Y(d_i))$ $(i\ins[1,r]$) generators of the ideal of $X\sst Y$. (Note that the $D_i$ are smooth on a neighborhood of $X$, considering the local embedding dimension.) So the product $\mprod_{i=1}^r\,H(d_i^{\one 2}z)$ gives the cohomology $L$-class of $N_{X/Y}$, see also \cite[Section 7]{LW} about a formula for the Pontryagin classes of smooth complete intersections. The last assertion follows from the Hirzebruch signature theorem, see Theorem\,\,\ref{T2.1}. (Recall that the degree of a $k$-dimensional subvariety of $\PP^N$ is defined to be the intersection number with a general $(N{-}k)$-dimensional linear subspace.) The last assertion follows from the equality
\begin{equation*}
H(d_i^{\one 2}z)^{-1}\eq\msum_{k=0}^m\,\bl(1\mi H(d_i^{\one 2}z)\br)^k\q\h{in}\,\,\,\Q[d,z]/(z^m),
\end{equation*}
since the constant term of $H(z)$ is 1. This finishes the proof of Proposition\,\,\ref{P2.1}.
\end{proof}

\begin{rem} \label{R2.3}
In the hypersurface case (that is, if $r\eq1$), we can verify for each concrete example that the last assertion of Proposition\,\,\ref{P2.1} is compatible with the Griffiths theorem on the Hodge numbers of projective hypersurfaces \cite{Gri}, which implies the information on signatures, see for instance \cite{JKSY0} for an explicit formula of Hodge numbers (where we may assume $X\sst\PP^{2m+1}$ is defined by $\msum_{i=1}^{2m+2}\,y_i^d\eq0$, since the Hodge numbers stay invariant by a deformation of smooth projective hypersurfaces). In the case $m\eq1$ with $d\eq4,5,6,7,8$, for instance, we see that the signatures are respectively 16, 35, 64, 105, 160 up to sign. These are the values of a numerical polynomial $-d(d\mi2)(d\pl2)/3$ as explained below. These assertions can be confirmed on a terminal of a computer on which Macaulay2 (preferably version 1.17 or later) is installed for instance as follows (after typing M2 and pressing Return):
\ms
\vbox{\fontsize{9pt}{4mm}\sf\pv#R=QQ[z]; d=4; m=1; T=QQ[t]; N=2*m+2; s=sub((t^(m+1)-1)/(t-1),T);#
\pv#ch=sub(1+z/2!+z^2/4!+z^3/6!+z^4/8!+z^5/10!+z^6/12!+z^7/14!+#
\pv#z^8/16!+z^9/18!+z^10/20!+z^11/22!+z^12/24!+z^13/26!+z^14/28!+#
\pv#z^15/30!+z^16/32!+z^17/34!+z^18/36!+z^19/38!+z^20/40!+z^21/42!,R);#
\pv#X=2*diff(z,ch); S=R/(z^(m+1)); H=sub(ch,S)*sub(s,{t=>1-sub(X,S)});#
\pv#G=sub(H,{z=>d^2*z})-1; P=H^N*sub(s,{t=>-G}); A=d*coefficient(z^m,P)#
\pv#T=QQ[t]; Q=sub(((t^d-t)/(t-1))^N,T); U=-t^d; R=sub((U^m-1)/(U-1),T);#
\pv#e=(m+1)*d; B=(-1)^m*coefficient(t^e,Q*(1+2*U*R))+1; A-B#}
\msn
It is quite surprising that the {\it order\one} of {\small\sf\verb#X=2*diff(z,ch);#} and {\small\sf\verb#S=R/(z^(m+1));#} {\it cannot\one} be reversed. Here the integer {\small\sf\verb#m#} must be {\it at most\one} 20 unless the definition of {\small\sf\verb#ch#} is expanded. One can also get $H^k$ mod ($z^{m+1}$) for $k$ at most 41 by typing {\small\sf\verb#H^k#} and pressing Return, where one has to set {\small\sf\verb#m=20#}. (One may have to add Return if it is erased by copying the above code.) The last part using the Griffiths theorem is valid for any positive integer {\small\sf\verb#m#} (although it may take very long if {\small\sf\verb#d#} is quite big). Here the sign {\small\sf\verb#(-1)^m#} is equal to $(-1)^{2m(2m-1)/2}$, and the last term {\small\verb#+1#} comes from the non-primitive part. In the first part, it is rather amazing that we can finally get an integer by calculating polynomials with coefficients whose denominators are quite large when {\small\sf\verb#m#} is not small.
\sk
We can replace {\small\sf\verb#R=QQ[z]; d=4;#} in the first line of the above code with {\small\sf\verb#R=QQ[d][z];#} where the last two lines must be deleted. We then get the numerical polynomials in variable {\small\sf\verb#d#} whose values give the signatures, where {\small\sf\verb#m#} must be at most 20 as before. (We say that $P(x)\ins\Q[x]$ is a {\it numerical polynomial\one} if $P(k)\ins\Z$ for any $k\ins\Z$, or equivalently, for any sufficiently large integer $k$. Here we take $Q(x)\eq aP(x)\ins\Z[x]$ with $a\ins\Z_{>0}$ and consider the values $Q(k)$ in the quotient ring $\Z/a\Z$.) Applying ``{\small\sf\verb#factor#}" to {\small\sf\verb#A#} for $m\eq1,2,3,4,5$ (for instance), the numerical polynomial {\small\sf\verb#A#} can be expressed respectively as
\begin{equation*}
\h{$\aligned-&\tfrac{1}{3}\one d(d\one ^2\mi4),\\&\tfrac{1}{15}\one d(2\one d\one ^4\mi10\one d\one ^2\pl23),\\-&\tfrac{1}{315}\one d(d\one ^2\mi4)(17d\one ^4\mi44\one d\one ^2\pl132),\\&\tfrac{1}{2835}\one d(62\one d\one ^8\mi510\one d\one ^6\pl1806\one d\one ^4\mi3590\one d\one ^2\pl5067),\\-&\tfrac{1}{155925}\one d(d\one ^2\mi4)(1382\one d\one ^8\mi8112\one d\one ^6\pl27018\one d\one ^4\mi41528\one d\one ^2\pl73215).\endaligned$}
\end{equation*}
These are generalized to the case $r\gess2$. For $r\eq2$, one has to replace {\small\sf\verb#QQ[d]#} by {\small\sf\verb#QQ[d1,d2]#}, and {\small\sf\verb#N#} must be {\small\sf\verb#2*m+3#}. Here one has to define {\small\sf\verb#G1#} and {\small\sf\verb#G2#}. One may propose the following.
\end{rem}

\begin{conj} \label{Con2.2}
The modified positivity $(-1)^jL^{2j}(X)\sgt0$ holds for any $j\ins[0,n/2]$ if $X$ is a smooth complete intersection for instance with $\msum_{i=1}^r\,d_i\gess n\pl r$ (at least).
\end{conj}
In the case of hypersurfaces with $n\eq2m$, one can verify that the signature multiplied by $(-1)^m$ is positive for $d\gess m$ if $m\ins[2,20]$ (and for $d\gess 3$ if $m\eq1$). This can be seen for instance by adding
\sk
{\small\sf\verb#(1/coefficient(d^(2*m+1),A))*A#}

\nin or

{\small\sf\verb#AA=sub(A,{d=>(2/3)*m*d}); (1/coefficient(d^(2*m+1),AA))*AA#}
\skn
and applying the following: A monic polynomial $P(d)\eq\msum_{k=0}^m\,(-1)^ka_kd\one ^{2m-2k+1}$ with $a_k\sgt 0$, $a_0\eq1$ has positive value at $d\in\Z_{\ges m}$ if for any odd $k\eq2j{+}1$ with $j\ins\bl[0,\tfrac{m{-}1}{2}\br]$, one has either $a_k\slt(m^2/2)^k $ or $a_k\slt m^2a_{k-1}$ with $k\gess3$. (Note that $(2/3)^2\slt 1/2$.) Here the first condition seems to be enough usually. (This would not be called a {\it mathematical\one} proof probably, since a computer is used.)

\begin{rem} \label{R2.4}
We can examine the above Conjectures\,\,\ref{Con2.1} (with $n\eq2m$) and \ref{Con2.2} for each case by using Singular \cite{Sing} as follows.
\sk
\vbox{\fontsize{9pt}{4mm}\sf\pv#LIB "general.lib"; ring R=0,(z),ds; int d,e,i,m; poly f,g,h,c,s; m=12; d=6; e=5;#
\pv#c=1; for(i=1; i<=m+1; i++){c=c+z^i/factorial(2*i);} s=2*diff(c,z); h=jet(jet(c,#
\pv#s,m),m); f=jet(h^(2*m+1),m); sprintf("L(P^
\pv#{if(subst(jet(f,i)-jet(f,i-1),z,1)<0){sprintf("Conj fails; (m,i)=(
\pv#exit;}} g=jet(jet(h^(2*m+2),subst(h,z,d^2*z),m),m); sprintf("L(X^{
\pv#2*m,d,g); for(i=1; i<=m; i++){if((-1)^i*subst(jet(g,i)-jet(g,i-1),z,1)<0)#
\pv#{sprintf("Conj fails; (m,d,i)=(
\pv#subst(h,z,d^2*z)*subst(h,z,e^2*z),m),m); sprintf("L(X^{
\pv#g); for(i=1; i<=m; i++){if((-1)^i*subst(jet(g,i)-jet(g,i-1),z,1)<0){sprintf#
\pv#("Conj fails; (m,d,e,i)=(
\skn
This code verifies Conjecture\,\,\ref{Con2.1} for $\PP^{2m}$ and Conjecture\,\,\ref{Con2.2} for a smooth hypersurface of degree $d$ and for a smooth complete intersection of bi-degree $(d,e)$ with $m,d,e$ given on the first line, where the dimension of the varieties are $2m$.
If the computation finishes without quitting, the conjectures are assured by the computer for these cases. (The usage of ``jet" to compute a quotient in the code is rather amazing, where one has to repeat ``jet" in the definition of $h,g$. One can compare the calculation of $L$-classes with the one by Macaulay2.)
Note that Conjectures~\ref{Con2.1} and \ref{Con2.2} are based	on computations using the above code.
\end{rem}

\section{L-class transformation} \label{S3}
Let $X$ be a compact variety with $A$ a subfield of $\R$. We say that a bounded complex with constructible cohomology sheaves $\Fb\ins D^b_c(X,A)$ is a {\it self-dual complex\one} if $\Fb$ is endowed with a perfect pairing
\begin{equation*}
S:\Fb{\otimes}_A\Fb\to\DD A_X,
\end{equation*}
where $\DD$ is the dual functor. Note that $\DD A_X\eq A_X[2d_X]$ if $X$ is smooth. Here the Tate twist is omitted choosing $\sqrt{-1}\ins\C$. We assume that $S$ is either {\it symmetric\one} or {\it skew-symmetric,} that is, $S\ssc\iota\eq S$ or $-S$, where $\iota$ is the involution $\Fb{\otimes}_A\G^{\ssb}\simto\G^{\ssb}{\otimes}_A\Fb$ defined by
\begin{equation}
\F^i{\otimes}_A\G^j\ni a{\otimes}b\mapsto(-1)^{ij}b{\otimes}a\in\G^j{\otimes}_A\F^i.
\label{3.1} \end{equation}
By $\Om_{A+}(X)$ and $\Om_{A-}(X)$ we mean respectively the {\it symmetric\one} and {\it skew-symmetric cobordism group\one} of $X$ with coefficients in $A$, see \cite{CS}, \cite{Ba1}. Note that we use {\it self-pairings\one} (rather than duality isomorphisms) as in \cite{FPS}. Put
\begin{equation*}
\Om_A(X)\eq\Om_{A+}(X)\oplus\Om_{A-}(X).
\end{equation*}
\sk
We define an $L$-{\it class transformation\one} using {\it cohomotopies\one} as in \cite{GM}, \cite{CS}, \cite{Ba1}. We may assume that cohomotopies $\rho\col X\tos S^{2k}$ are differentiable maps by the smooth approximation theorem as in \cite{Ba1} if $X$ is smooth. In general we can define differentiable functions on complex varieties using minimal local embeddings into smooth varieties.
\sk
Let $s\ins S^{2k}$ be a general point so that the restriction of $\rho$ to any stratum $Z$ of a Whitney stratification of $X$ is a submersion on a neighborhood of $X_s\cap Z$, where $X_s\defs\rho^{-1}(q)$. Set $\Fsb\defs\Fb[-k]|_{X_s}$. Since $(\DD A_X)|_{X_s}\eq(\DD A_{X_s})[2k]$, we get the induced pairing
\begin{equation}
S|_{X_s}:\Fsb{\otimes}_A\Fsb\to\DD A_{X_s},
\label{3.2} \end{equation}
by using the canonical isomorphisms
\begin{equation}
\aligned(&\Fb{\otimes}\Fb)|_{X_s}=A_{X_s}[k]{\otimes}\Fsb{\otimes}A_{X_s}[k]{\otimes}\Fsb\\&=A_{X_s}[k]{\otimes}A_{X_s}[k]{\otimes}\Fsb{\otimes}\Fsb\\&=A_{X_s}[1]{\otimes}A_{X_s}[1]{\otimes}\cdots{\otimes}A_{X_s}[1]{\otimes}A_{X_s}[1]{\otimes}\Fsb{\otimes}\Fsb.\endaligned
\label{3.3} \end{equation}
Here $\Fsb[k]$ is identified with $A_{X_s}[k]{\otimes}_A\Fsb$, the second and third isomorphism follow from \eqref{3.1}, and the third one produces the sign $(-1)^{k(k-1)/2}$, see \cite[(1.1.4.2)]{De1}, \cite{De2}. The last assertion may be compared with the identity
\begin{equation*}
\omega_1{\wedge}\cdots{\wedge}\one\omega_k{\wedge}\one\overline{\omega}_1{\wedge}\cdots{\wedge}\one\overline{\omega}_k=(-1)^{k(k-1)/2}\omega_1{\wedge}\one\overline{\omega}_1{\wedge}\cdots{\wedge}\one\omega_k{\wedge}\one\overline{\omega}_k,
\end{equation*}
with $\omega_1,\dots,\omega_k$ one-forms, which is closely related to the positivity of a polarization on the primitive part especially in the abelian variety case.
\sk
We define a modified $L$-class $\Lt_k(\Fb,S)\ins H_{2k}(X,\Q)$ using the {\it signature\one} of the symmetric bilinear form
\begin{equation}
S|_{X_s}:H^0(X_s,\Fsb){\otimes}_AH^0(X_s,\Fsb)\to A,
\label{3.4} \end{equation}
for each cohomotopy $\rho\col X\to S^{2k}$, where $k$ is even if $S$ is symmetric, and odd otherwise (and $\Lt_k(\Fb,S)\eq0$ in the other cases). Note that $X$ must be replaced by its product with $S^m$ for $m$ sufficiently large in the case $4k\mi1\less 2d_X$, see \cite{Ba1}. For the proof of the independence of the choices of cohomotopy and $s\ins S^{2k}$, one can use the base change of the diagrams
\begin{equation*}
\xymatrix{A_I \ar[r] \ar[d] & A_{\{v\}} \ar[d] & A_{S^{2k}} \ar[r] \ar[d] & A_{\{s'\}} \ar[d]\\ A_{\{u\}} \ar[r] & j_!A_{I^{\circ}}[1] & A_{\{s\}} \ar[r] & j'_!A_{S^{2k}\setminus\{s,s'\}}[2k]}
\end{equation*}
Here $I^{\circ}\defs(u,v)\buildrel{j}\over{\into}I\defs[u,v]\sst\R$ and $j'\col S^{2k}\setminus\{s,s'\}\into S^{2k}$ are canonical inclusions with $s,s'$ general points of $S^{2k}$ and cohomotopy fixed, and $I$ is a {\it very small\one} closed subinterval of the parameter space of cohomotopies with general point $s\in S^{2k}$ fixed. (The set of regular values of a differentiable cohomotopy is not necessarily connected.) Note that the upper-left terms of the diagrams are dual of the lower-right terms, where the direct images under closed embeddings are omitted, see \cite[(1.1.4-5)]{FPS}.
\sk
In the case $X$ is irreducible, we define
\begin{equation*}
\Lt_{\ssb}(X)\defs\Lt_{\ssb}({\rm IC}_XA,S_X)\,\,\,\,\h{in}\,\,\,H_{2\ssb}(X,\Q),
\end{equation*}
where $S_X$ is the natural self-pairing of the intersection complex ${\rm IC}_XA$ (containing no sign). Note that
\begin{equation}
\Lt_{d_X-k}(X)\eq0\,\,\,\h{unless}\,\,\,k\,\,\h{is even,}
\label{3.5} \end{equation}

\begin{rem} \label{R3.1}
For $k\in\N$ we have the equalities
\begin{equation}
(-1)^{k(k-1)/2}\Lt_k(X)\eq L_k(X)\,\,\,\,\h{in}\,\,\,H_{2k}(X,\Q).
\label{3.6} \end{equation}
Here it is enough to consider cohomotopies $\rho\col X\tos S^{2k}$ only in the case $d_X{-}\one k$ is {\it even.} So a sign does {\it not\one} appear in the second isomorphism of \eqref{3.3} if an intersection complex underlies a pure Hodge module associated with a constant variation of Hodge structure with fiber-wise weight {\it even.} Indeed, we may assume that $s\ins S^{2k}$ is sufficiently general, and the fiber $\rho^{-1}(s)$ intersects transversally every stratum of a Whitney stratification of $X$. The sign is then determined by the restriction to the smooth part, and the (real) dimension $d_{X_s}$ of the stratified space $X_s$ is {\it divisible by\one} 4 with $\Fsb\eq A_{X_s}[d_{X_s}/2]$ in the smooth case. In the case an intersection complex underlies a pure Hodge module associated with a variation of Hodge structure with fiber-wise weight {\it odd,} however a sign {\it does\one} appear when the complex dimension of the support is even and $k$ is odd, see also \cite[Example 3.22]{AMS} (where duality isomorphisms are mainly used). This can be applied to direct factors of direct images of constant sheaves by projective morphisms, and it seems rather difficult to get the functoriality {\it without any sign.} The sign convention in \cite{CS}, \cite{BSY} does not seem very clear to non-specialists, since duality isomorphisms are used instead of self-pairings, see \cite{FPS}. In this paper we use the $L$-class transformation given by $(-1)^{k(k-1)/2}\Lt_k$, which is compatible with direct images (see Remark~\ref{R3.2} below) and coincides with the Goresky-MacPherson $L$-class, see also \cite{AMS} for more details.
\end{rem}

\begin{rem} \label{R3.2}
The above $L$-class transformation is compatible with the pushforward by a morphism of compact varieties $f\col X\tos Y$, that is, for $[(\Fb,S)]\ins\Om_{A,\pm}(X)$, we have
\begin{equation}
f_*\Lt_{\ssb}(\Fb,S)=\Lt_{\ssb}\bl(f_*(\Fb,S)\br)\,\,\,\,\h{in}\,\,\,H_{2\ssb}(Y,\Q).
\label{3.7} \end{equation}
Indeed, let $\rho'\col Y\to S^{2k}$ be a cohomotopy which is transversal to any stratum of a Whitney stratification of $Y$ compatible with one of $X$ via $f$. Let $s\ins S^{2k}$ be a  general point. Set $\rho\defs\rho'\ssc f$. Replacing $\Fb$ with the Godement canonical flasque resolution, we may assume that the $\F^j$ are flasque sheaves and vanish for $j\ll0$. Define $X_s$, $\Fsb$ and $Y_s$ as above. Consider the isomorphism obtained by replacing $\Fb$, $X_s$ and $\Fsb$ with $f_*\Fb$, $Y_s$ and $f_*\Fsb$ respectively in \eqref{3.3}. We have a canonical morphism from this isomorphism to the direct image of \eqref{3.3} producing a commutative diagram. The assertion then follows from this commutativity considering also the direct image of the canonical flasque resolution of \eqref{3.3}.
\end{rem}

\begin{rem} \label{R3.3}
By \eqref{3.7} the compatibility of the $L$-class transformation with the pushforward under a proper morphism would not hold without giving the sign properly, for instance in the case $X$ is smooth, $\dim X$ is odd, and $\Hc^0\!f_*A_{h,X}[d_X]$ contains a Hodge module such that the dimension of its support is even and non-zero, where some sign should appear for the construction of the $L$-class as explained at the end of Remark\,\,\ref{R3.1}. Note that the sign does not appear in the construction on $X$ by \eqref{3.6}. It does not seem that the two constructions on $Y$ using the sign and without using the sign coincide in general.
\sk
Set for instance $W\defs\{s(x^3{+}y^3{+}z^3)\eq txyz\}\sst\PP^2{\times}S$ with $S\defs\PP^1$. Put $Y\defs Y'{\times}D$ with $Y'$ the relative projective cone of $Z\defs W{\times}C$ over $S$, where $C\defs\PP^1$ and $D$ is a compact smooth variety with $\dim_{\C}D$ odd. Set $X\defs X'{\times}D$ with $X'$ the blowup along the vertex $S$ of the cone $Y'$ (which is a $\PP^1$-bundle over $Z$). Here we consider the composition of the projection $Y\onto D$ with a cohomotopy $D\onto S^{2k}$, where $k\eq\dim_{\C}D$. We need the product with $C\,({\cong}\,\PP^1)$, since the ``relative non-primitive part" gives the direct factor of the direct image by the blowup which is supported on $S$, using the Thom-Gysin sequence and shrinking $S$. There are four singular fibers of the one-parameter family $W\onto S$, which are unions of three lines having three singular points. Using the Leray-type spectral sequence and the classical Hodge index theorem for surfaces, it is enough to calculate the Euler number of $W$ (which is given by $4\,{\cdot}\,3\eq12$) in order to show the non-vanishing of the signature, since $W$ is birational to $\PP^2$ (putting $z\eq s\eq1$) so that $H^0(W,\Om_W^p)\eq0$ for $p\eq1,2$. Here a constant variation does not work, since we get the signature 0.
\end{rem}

\begin{prop} \label{P3.4}
In the notation of {\rm\eqref{1.1}}, set $L^{\rm HT}_*(X)\defs T_{1*}({\rm IC}_X\R_h[-d_X])$, called the Hodge-theoretic intersection $L$-class. We have the equality $L_{\ssb}^{\rm HT}(X)\eq L_{\ssb}(X)$ in $H_{2\ssb}(X,\Q)$ if $\dX\less1$.
\end{prop}

\begin{proof}
The image of $sd([X\,{\buildrel{\rm\!\!id}\over\to}\,X])$ by the topological $L$-class transformation can be calculated by applying $(-1)^{k(k-1)/2}\Lt_k$ to the cubic hyperresolution of $X$, and coincides with $T_{1*}(X)$ by Remarks\,\,\ref{R3.1} and \ref{R3.2} using \cite{Hi} (see also Remark\,\,\ref{R2.1}). The hypothesis then implies that the difference between $sd([X\,{\buildrel{\rm\!\!id}\over\to}\,X])$ and $[({\rm IC}_X\R,S)]$ is given by a $\Z$-linear combination of the classes of self-dual complexes $(\Fb_i,S_i)$ which underlie polarized pure $\R$-Hodge modules $\M_i$ whose support $Z_i$ has dimension at most 1, where $S_i$ is a polarization of Hodge module. (This coincides with a polarization of variation of Hodge structure, since $\dim Z_i\eq1$ or $0$, see \cite[Section 5.4.1]{mhp}.)
\sk
It is then enough to show the equalities
\begin{equation}
T_{1,k}(\M_i)=\Lt_k(\Fb_i,S_i)\,\,\,\,\h{in}\,\,\,H_{2k}(X,\Q),
\label{3.8} \end{equation}
for $k\eq0,1$, where $(-1)^{k(k-1)/2}\eq1$.
If $k\eq0$, the assertion follows from \cite[Theorem 5.3.1]{mhp} applied to the direct image by $a_{Z_i}\col Z_i\tos pt$. (Note that $H^0(a_{Z_i})_*$ preserves the weight.)
For $k\eq1$, we have to calculate the sum of the local intersection numbers of $Z_i$ with a general fiber of a cohomotopy $\rho$ which are multiplied by the signature of the Hodge structure at each intersection point (which is independent of the intersection point since $Z_i$ is irreducible) where $\dim Z_i\eq1$. We may assume that the intersection is transversal, since $s\ins S^2$ is general. The local intersection numbers are $\pm1$, and given by choosing an orientation of $S^2$. Note that this orientation is employed also in the definition of the Hurewicz map.
\sk
We will denote $Z_i$ by $Z$ to simplify the notation. Let $\Zt\tos Z$ be the normalization. Let $\rht_Z\col\Zt\tos S^2$ be its composition with the restriction of $\rho$ to $Z$. Let $\rht_{Z_s}$ be the restriction of $\rht_Z$ to $\Zt_s\defs\rht_Z^{-1}(s)$ for a general point $s\ins S^2$. We have the commutative diagram
\begin{equation} \begin{gathered} \xymatrix{
(\rht_Z)_*A_{\Zt} \ar[r] \ar[d]^{\rm Tr} & (\rht_{Z_s})_*A_{\Zt_s} \ar[d] \\
A_{S^2} \ar[r] & A_{\{s\}} } \end{gathered}
\label{3.9} \end{equation}
where the left vertical morphism is the {\it trace\one} morphism, which is defined by using Verdier duality, and the right one is its base change over $s$. Taking the global sections, we see that the right-hand side of \eqref{3.8} is given by the homology class of $Z_i$ multiplied by the signature of a polarization (which coincides with the Hodge signature) of a general stalk of the generic variation of Hodge structure of $\M_i$. (Note that $(-1)^{d_X(d_X-1)/2}=(-1)^n$ if $d_X\eq2n$, see also a remark after \cite[Theorem 3]{FPS}.) This is the same for the left-hand side of \eqref{3.8} applying \eqref{4.1} in Section\,\,\ref{S4} below to a desingularization of $Z_i$. So the assertion follows. This finishes the proof of Proposition\,\,\ref{P3.4}.
\end{proof}

The following is pointed out in \cite{AMS} (where arguments using duality isomorphisms rather than self-pairings are explained).

\begin{prop} \label{P3.5}
The degree $2\dX$-part of the $L$-classes $L_{\dX}(X)$ and $L^{\rm HT}_{\dX}(X)$ coincide.
\end{prop}

\begin{proof}
By an argument similar to the beginning of the proof of Proposition\,\,\ref{P3.4}, the assertion can be reduced to the following:

\begin{claim}
For a pure $\R$-Hodge module $\M$ with underlying self-dual $\R$-complex $(\Fb,S)$, an equality similar to \eqref{3.8} holds for $k\eq d$ (with $d\defs\dim{\rm Supp}\,\M$) by adding the sign $(-1)^{d(d-1)/2}$ to the right-hand side.
\end{claim}

We may assume ${\rm Supp}\,\M$ is smooth using the compatibility of the class transformations with pushforwards under proper morphisms. In the case $\dim X$ is odd, we can verify that both sides $T_{1,d}(\M)$ and $(-1)^{d(d-1)/2}\Lt_d(\Fb,S)$ vanish using arguments similar to the proofs of Proposition\,\,\ref{P3.4} and Theorem\,\,\ref{T1} in Section\,\,\ref{S7} below. It is relatively easy to see that these coincide up to sign by a similar argument in the case $\dim X$ is even. We can examine the sign using the trivial identity
\begin{equation*}
(n\pl d)(n\pl d\pl 1)/2-d(d\mi 1)/2\,\equiv\,n/2\pl d\mod2.
\end{equation*}
Here the left-hand side comes from \cite[Theorem 3]{FPS} and \cite[Section 5.4.1]{mhp}, and calculates the sign of $(-1)^{d(d-1)/2}\Lt_d(\Fb,S)$, which employs the ``naive" restriction to points (not using a sign), and coincides with the signature of the polarization of Hodge structure at the points ``up to sign", and this sign is calculated by the left-hand side. (Note that spheres are not naturally oriented, but the obtained sign is independent of the choice, since the Hurewicz map depends on it as well as the induced orientation on the fiber of a cohomotopy.)
\sk
We see that $n/2$ on the right-hand side corresponds to the Tate twist $(n/2)$, which can be applied to get a variation of Hodge structure of fiber-wise weight 0, and induces the corresponding shift of the Hodge filtration, see also Lemma\,\,\ref{L3.6} below. Moreover $d$ corresponds to the shift of the de Rham complex in the $\D$-module theory as remarked in the proof of Theorem\,\,\ref{T1} in Section\,\,\ref{S7}. The assertion is then reduced to the following (which seems implicit in \cite{Sc}).
\end{proof}

\begin{lem} \label{L3.6}
Let $(H_{\R},(H_{\C},F),S)$ be a polarized $\R$-Hodge structure of type $\{(p,q),(q,p)\}$ such that $\dim_{\R}H_{\R}\eq2$ and $n\defs p\pl q$ is even with $p\sgt q$. Then the signature of the symmetric $\R$-bilinear form $(H_{\R},S)$ coincides with $(-1)^{p-n/2}\one2$ $($which is invariant by Tate twists$)$.
\end{lem}

\begin{proof}
Take a non-zero $v\ins F^pH_{\C}$. By the definition of polarization, we have
\begin{equation*}
i^{q-p}S_{\C}(v,\overline{v})\eq(-1)^{p-n/2}S_{\C}(v,\overline{v})\sgt0,\q S_{\C}(v,v)\eq0,
\end{equation*}
where $S_{\C}$ is the scalar extension of $S$ (and $i^{q-p}\eq i^{p-q}$ since $p\pl q$ is even).
\sk
Setting $u_1\defs(v\pl\overline{v})/2$, $u_2\defs (v\mi\overline{v})/2i\in H_{\R}$, we verify that
\begin{equation*}
(-1)^{p-n/2}S(u_j,u_j)\sgt0\q\h{for}\,\,\,j=1,2,\q\h{and}\q S(u_1,u_2)\eq0.
\end{equation*}
So the assertion follows.
\end{proof}

\section{Topological filtration} \label{S4}
The Grothendieck group $K_0(X)$ of coherent sheaves on a complex algebraic variety $X$ has the {\it topological filtration\one} (denoted by $G$ in this paper) such that the $G_kK_0(X)$ are generated by the classes of coherent sheaves $\G$ with $\dim{\rm Supp}\,\G\less k$ (see \cite[Examples~1.6.5 and 15.1.5]{Fu}, \cite{SGA6}), and the Todd class transformation $td_*$ to the Chow groups induces the isomorphisms
\begin{equation}
\aligned td_*:K_0(X)_{\Q}&\simto\mopl_k\,{\rm CH}_k(X)_{\Q},\\
\Gr^G_k\one td_*:\Gr_k^GK_0(X)_{\Q}&\simto{\rm CH}_k(X)_{\Q},\endaligned
\label{4.1} \end{equation}
see \cite[Corollary 18.3.2]{Fu}. (We may assume $X$ is smooth for the applications in this paper.) The increasing filtration $G_j$ on $\mopl_k\,{\rm CH}_k(X)_{\Q}$ is defined by taking the direct sum over $k\less j$, and the inverse of the last isomorphism of \eqref{4.1} is given by $[Z]\mapsto[\OO_Z]$ for irreducible reduced closed subvarieties $Z$ of $X$ with dimension $k$. The Todd class transformation $td_*$ to the Borel-Moore homology is identified with the {\it cycle class map\one} taking the graded quotients of the topological filtration $G$.
Note that passing to the graded pieces of $G$ is quite different from restricting to some open subvariety; the latter loses a lot of information in general except for the top degree. (In \cite[Lemma 4.16]{AMS} one could assume $X$ is sufficiently small so that the vector bundles in the lemma are trivial and the argument is quite simplified.)
\sk
If $X\eq\PP^m$, for a coherent sheaf $\G$ on $\PP^m$, we have
\begin{equation}
\Gr_k^G[\G]=\msum_i\,m_i\deg Z_i\q\h{in}\,\,\,\,\Gr_k^GK_0(\PP^m)_{\Q}\eq{\rm CH}_k(\PP^m)_{\Q}\eq\Q.
\label{4.2} \end{equation}
Here $k\eq\dim{\rm supp}\,\G$, and the $Z_i$ are maximal-dimensional irreducible components of the support of $\G$ with $m_i\ins\Z_{>0}$ the multiplicity of $\G$ at the generic point of $Z_i$. Note that the degree of $Z_i$ is defined to be the intersection number of $Z_i$ with a general linear subspace of the complementary dimension, see also \cite[1.3]{MSS2}, \cite[2.2]{MSY}.

\section{Non-characteristic restrictions} \label{S5}
For a complex variety $X$, consider the mapping cone
\begin{equation*}
\M_X^{\ssb}:=C\bl(\R_{h,X}\to\mopl_i\,{\rm IC}_{X_i}\R_h[-d_{X_i}]\br)\in D^b\MHM(X,\R),
\end{equation*}
where the $X_i$ are irreducible components of $X$. Note that
\begin{equation}
{\rm Supp}\,\M_X^{\ssb}\defs\mcup_j\,{\rm Supp}\,H^j\M^{\ssb}\eq\Si_X,
\label{5.1} \end{equation}
where $H^{\ssb}$ denotes the standard cohomology functor $D^b\A\tos\A$ for $\A\eq\MHM(X)$. Indeed, it is quite well-known that $X$ is a $\Q$-homology manifold if and only if $\M_X^{\ssb}\eq0$. (This equivalence is an easy consequence of \cite[(4.5.6--9)]{mhm}.)
\sk
Let $X'$ be a {\it stratified-smooth\one} hypersurface of $X$ as in Proposition\,\,\ref{P1} with $i_{X'}\col X'\into X$ the inclusion. Set $X'_i\defs X_i\cap X'$. We have the following.
\begin{lem} \label{L5.1}
There are isomorphisms in $D^b(X',\R)\,{:}$
\begin{equation}
i_{X'}^*{\rm IC}_{X_i}\R={\rm IC}_{X'_i}\R[1],\q i_{X'}^*\DD(\R_X)=\DD(\R_{X'})(1)[2],
\label{5.2} \end{equation}
and an isomorphism in $D^b\MHM(X')\,{:}$
\begin{equation}
i_{X'}^*\M_X^{\ssb}=\M_{X'}^{\ssb}:=C\bl(\R_{h,X'}\to\mopl_i\,{\rm IC}_{X'_i}\R_h[-d_{X'_i}]\br).
\label{5.3} \end{equation}
\end{lem}

(Here $i_{X'}^*$ denotes also the pullback in the bounded derived categories of mixed Hodge modules, and $\DD$ denotes the dual functor.)

\begin{proof}
The isomorphisms in \eqref{5.2} follow from the transversality of $X'$ to the strata of a Whitney stratification. The first isomorphism of \eqref{5.2} implies \eqref{5.3}. Indeed, the intersection complex Hodge module ${\rm IC}_{X'_i}\R_h$ is characterized up to Tate twist by the condition that its underlying $\R$-complex is the intersection complex ${\rm IC}_{X'_i}\R$. Here we use the functorial morphisms $j_!j^*\tos id\to j_*j^*$ factorizing the morphism in the definition of intermediate direct image \cite{BBD} with $j$ the inclusion of a non-empty smooth affine open subvariety of $X'_i$. (Note that the above condition implies that it is a {\it simple\one} Hodge module.) This finishes the proof of Lemma\,\,\ref{L5.1}.
\end{proof}
\sk
Let $Y$ be a smooth affine variety containing $X$ locally with $Y'\sst Y$ a smooth hypersurface defined by a function $z$ on $Y$ such that $X'\cap Y\eq X\cap Y'$. Let $i_{Y'}\col Y'\into Y$ be the natural inclusion. Let $(M_{X,Y}^{\ssb},F)$ be the underlying complex of filtered holonomic left $\D_Y$-module of $\M^{\one\ssb}{}\!\!\!_X|_Y$, and similarly for $(M_{X',Y'}^{\ssb},F)$.

\begin{prop} \label{P5.1}
With the above notation, the inclusion $i_{Y'}$ is non-characteristic to the filtered cohomology $\D_Y$-modules $\Hc^j(M_{X,Y}^{\ssb},F)$. More precisely, the $V$-filtration along $Y'$ on the $\Hc^j\!M_{X,Y}^{\ssb}$ is given by the $z$-adic filtration, the action of $z$ on the $\Gr^F_p\Hc^j\!M_{X,Y}^{\ssb}$ is injective for $p,j\ins\Z$, and we have the isomorphisms of filtered holonomic left $\D_{Y'}$-modules
\begin{equation}
\Hc^j(M_{X',Y'}^{\ssb},F)=\OO_{Y'}{\otimes}_{\OO_Y}\Hc^{j+1}(M_{X,Y}^{\ssb},F)\q(j\ins\Z).
\label{5.4} \end{equation}
\end{prop}

\begin{proof}
The pullback by the closed immersion $i_{Y'}$ can be calculated by the mapping cone
\begin{equation*}
C({\rm can}:\psi_{z,1}\to\varphi_{z,1}),
\end{equation*}
see \cite[(2.24.2)]{mhm}. This definition of $i_{Y'}^*$ is compatible with taking the cohomology functor $\Hc^{j+1}$ on the right-hand side of \eqref{5.4}, since the nearby and vanishing cycle functors are {\it exact\one} functors of mixed Hodge modules. Moreover the vanishing cycle functor $\varphi_z$ applied to the $\Hc^j\!M_{X,Y}^{\ssb}$ vanish by the assumption on the $Y'$. So the $V$-filtration on $\Hc^j\!M_{X,Y}^{\ssb}$ along $z$ is given by the $z$-adic filtration, and the $\Gr_V^{\al}$ vanish unless $\al\ins\Z_{>0}$. The Hodge filtration $F$ is not shifted when we apply $\psi_z$, since we do not take the graph embedding for $z$. The assertion about the injectivity of the action of $z$ is then verified by using \cite[(3.2.1.2)]{mhp} with the nine lemma, see also \cite{DMST}. This finishes the proof of Proposition\,\,\ref{P5.1}.
\end{proof}

\section{Proof of Proposition~\ref{P1}} \label{S6}
Using the last isomorphism of \eqref{5.2} in Lemma\,\,\ref{L5.1} together with the functorial morphism $id\tos(i_{X'})_*i_{X'}^*$, we get the canonical morphism
\begin{equation}\label{6.1}
H_j(X,\Q)\to H_{j-2}(X',\Q)(1).
\end{equation}
It is easy to see its compatibility with the cycle class map for cycles which are transversal to $X'$ (using a Whitney stratification of cycles), see \cite{Fu} for the general case.
\sk
We may assume $\dX\gess1$, since otherwise the assertion is trivial.
Let $U$ be an smooth open subvariety of $\Si_X$ whose closure $\overline{U}$ is the union of $\dX$-dimensional irreducible components of $\Si_X$. Take a smooth compactification $Z$ of $U$ having a projective morphism $\pi\col Z\tos X$ extending the locally closed embedding $i_U\col U\into X$.
\sk
Let $(M_U^{\ssb},F)$ be the underlying complex of filtered holonomic $\D_U$-module of $i_U^*\M^{\ssb}{}\!\!\!_X$. The latter complex can be extended to a bounded complex of mixed Hodge modules $\M_Z^{\ssb}$ on $Z$ (using the open direct image, for instance). Let $(M_Z^{\ssb},F)$ be the underlying complex of filtered holonomic $\D_Z$-modules of $\M_Z^{\ssb}$. Let $G$ be the topological filtration in Section\,\,\ref{S4}. We have the following equalities in $\Gr^G_{\dX}K_0(X)[y,y^{-1}]$\,:
\begin{equation}\label{6.2} 
\aligned&\Gr^G_{\dX}\DR_y[\M_X^{\ssb}]=\pi_*\Gr^G_{\dX}\DR_y[(M_Z^{\ssb},F)]\\={}&\msum_{j\in\Z}\,(-1)^j\pi_*\Gr^G_{\dX}\DR_y[(\Hc^j\!M_Z^{\ssb},F)]\\={}&\msum_{j,p\in\Z}\,(-1)^{j+\dX}\pi_*\Gr^G_{\dX}[\Gr^F_{-p}\Hc^j\!M_Z^{\ssb}](-y)^p(1\pl y)^{\dX}.\endaligned
\end{equation}
Indeed, $(1\pl y)^{\dX}\eq\msum_{j=0}^{\dX}\,\binom{\dX}{j}y^j$ with $\binom{\dX}{j}\eq{\rm rk}\,\Om_Z^j$ (so its multiplication corresponds to the passage from filtered $\D_Z$-modules to the associated de~Rham complexes on $Z$), where the differential of a {\it bounded\one} complex can be {\it neglected\one} (that is, it can be replaced by 0) in the Grothendieck group. We define $\DR_y[(M_Z^{\ssb},F)]$ in the same way as \eqref{1.1}, and $\pi_*\DR_y[(M_Z^{\ssb},F)]$ may be replaced by any coherent extension of $\DR_y[(M_U^{\ssb},F)]$. Note that the additional sign $(-1)^{\dX}$ in the last term of \eqref{6.2} comes from the shift of the de Rham complex in $\D$-module theory (more precisely, $\Hc^k\DR_U(\Hc^jM_U^{\ssb})\eq0$ for $k\nes{-}\dX$, $j\ins\Z$).
\sk
The assertion then follows from Propositions\,\,\ref{P3.5} and \ref{P5.1} taking the graded piece $\Gr^G_{\dX}$ of the Todd class transformation, since we get the cycle class map as is noted after \eqref{4.1}. This finishes the proof of Proposition\,\,\ref{P1}. \hfill\qedsymbol

\begin{rem} \label{R6.1}
It is not difficult to construct a counterexample to the converse of the first assertion of Proposition\,\,\ref{P1}. Consider for instance a hypersurface $X\sst Y\defs\PP^4{\times}\PP^1$ such that
\begin{equation}\label{6.3}
{\rm Sing}\,X\eq\{0\}{\times}\PP^1\cup C{\times}\{\infty\},
\end{equation}
and the transversal slice to the second component is type $A_1$, while the one to the first is a $\mu$-constant deformation of a singularity defined by
\begin{equation*}
g\defs x'{}^7\pl y'{}^5\pl z'{}^4\pl v'{}^3\pl x'y'z'v'.
\end{equation*}
Here $C\,({\cong}\,\PP^1)\sst\PP^4$ is an intersection of three hyperplanes of $\PP^4$, and $(x',y',z',v')$ is a local coordinate system. (Note that $\tfrac{1}{7}\pl\tfrac{1}{5}\pl\tfrac{1}{4}\pl\tfrac{1}{3}\slt1$.) Its defining polynomial can be written as
\begin{equation*}
\aligned h:={}&(\al x^7\pl\beta y^7\pl\ga z^7\pl\delta v^7\pl y^5w^2\pl z^4w^3\pl v^3w^4\pl xyzvw^3)s^2\\
&+\bl(\al'x^7\pl\beta'y^7\pl\ga'z^7\pl(x^2\pl y^2\pl z^2)v^5\br)t^2,\endaligned
\end{equation*}
with $\al,\beta,\ga,\delta,\al',\beta',\ga'$ sufficiently general complex numbers. We see that the restrictions of $X$ to the open subsets
\begin{equation*}
\{xs\ne0\},\,\{ys\ne0\},\,\{zs\ne0\},\,\{vs\ne0\},\,\{xt\ne0\},\,\{yt\ne0\},\,\{zt\ne0\}
\end{equation*}
are smooth, substituting $x\eq s\eq1$, etc., since a general fiber of a morphism from a smooth variety is smooth (and the intersection of a finite number of dense Zariski-open subsets is non-empty). This implies \eqref{6.3}, where $C\eq\{x\eq y\eq z\eq0\}\sst\PP^4$.
\sk
Substituting $w\eq s\eq1$ and $v\eq t\eq1$, we get respectively
\begin{equation*}
\aligned&\al x^7\pl y^5\pl z^4\pl v^3\pl xyzv\pl\h{higher terms},\\&x^2\pl y^2\pl z^2\pl \delta s^2\pl\h{higher terms},\endaligned
\end{equation*}
where the ``higher terms" are with respect to the Newton filtrations.
\sk
The unipotent monodromy part of the vanishing cohomology of $g$ is three-dimensional and has only one Jordan block. This can be verified for instance by using \cite{JKSY} or Singular \cite{Sing}. For the $A_1$-singularity, it is one-dimensional. The Hodge signatures of their link cohomologies are then $\pm 1$ with {\it opposite\one} signs.
\sk
The difference of the two $L$-classes on $X$ at degree 2 is given by a linear combination of the cycle classes of $\{0\}{\times}\PP^1$ and $C{\times}\{\infty\}$ (where the coefficients are $\pm 1$ with {\it opposite\one} signs using Theorem\,\,\ref{T1}). They are linearly independent, since their images in $H_2(Y,\Q)$ are. We see that a {\it cancellation\one} occurs by restricting to a general hyperplane section $X'$ which is rationally equivalent to the union of the pullbacks of hyperplanes of $\PP^4$ and $\PP^1$ (using the Segre embedding). So this gives a counterexample.
\end{rem}

\section{Proof of Theorem~\ref{T1}} \label{S7}
We may assume $\dX\gess1$, since the case $\dX\eq0$ follows from the argument in the proof of \cite[Theorem 2]{FPS} (that is, Theorem\,\,\ref{T2} in this paper).
Using the equalities in \eqref{6.2} as in the proof Proposition~\ref{P1}, Theorem\,\,\ref{T1} follows from Propositions\,\,\ref{P3.5} and \ref{P5.1} together with \eqref{4.1} (applied to $Z$) and also the last remark in Section\,\,\ref{S4}.
\sk
Note that the modified signature appears from the definition of ${\rm DR}_y$ after the restriction to a point of $U$ (see also a remark after Theorem\,\,\ref{T2}), and the filtered complex $(M_{X,Y}^{\ssb},F)$ is isomorphic to the direct image of $(M_U^{\ssb},F)$ as a filtered $\D$-module by the closed immersion $U\into Y$, where we may assume $U\eq\overline{U}\cap Y$ (replacing $U,Y$ if necessary). This finishes the proof of Theorem\,\,\ref{T1}.\hfill\qedsymbol

\section{Proof of Theorem~\ref{T3}} \label{S8}
Since $Z$ is sufficiently general and $X$ is very ample, we may assume that $Z$ intersects $X$ transversally and $Z\cap{\rm Sing}\,X\eq\emptyset$. Choosing a smooth path $\ga$ connecting $0$ and $1$ in the complement of the discriminant of $g$ in $\PP^1$ (except the origin), where $g^{-1}(0)\eq X$, $g^{-1}(1)\eq Z$, we get morphisms $r^*\col H^{\ssb}(X)\tos H^{\ssb}(Z)$, $r_*\col H_{\ssb}(Z)\tos H_{\ssb}(X)$ which are induced by a retraction map $r\col Z\to X$, see Remark~\ref{R8.1} below. Here one can use also the composition
\begin{equation*}
r^*:H^{\ssb}(X)\eq H^{\ssb}\bl(g^{-1}(\ga)\br)\tos H^{\ssb}(Z),
\end{equation*}
and similarly for $r_*$ using duality.
\sk
Set $V\defs Y{\times}\PP^1$. We have a natural inclusion $\Yt\sst V$. Put $X_c\defs\Yt\cap Y_c$ with $Y_c\defs Y{\times}\{c\}$ for $c\ins\PP^1$, where $X_0\eq X$, $X_1\eq Z$. We have $[T\Yt]\eq[TV|_{\Yt}]\mi[N_{\Yt/V}]$, hence we get the equality
\begin{equation*}
[T_{\rm vir}X_c]\defs[TY_c|_{X_c}]\mi[N_{X_c/Y_c}]\eq[T\Yt|_{X_c}]\mi[N_{X_c/\Yt}]\q\h{in}\,\,\,K^0(X_c).
\end{equation*}
Indeed, $[TV|_{Y_c}]\eq[TY_c]\pl[N_{Y_c/V}]$ and $N_{Y_c/V}|_{X_c}\eq N_{X_c/\Yt}$.
Note also that $\Yt$ is smooth, since $\Yt$ is the blowup of $Y$ along $X\cap Z$ with $Z$ sufficiently general and $X$ has only isolated singularities.
\sk
Set $\M\defs\psi_{g,1}\R_{h,\Yt}[d_X]$. The weight filtration $W$ on $\M$, $V_{x,1}$, $V_{x,\ne1}$ is the monodromy filtration shifted by $d_X$, $d_Y\eq d_X{+}1$, and $d_X$ respectively, and we have the isomorphisms for $k\sgt0$\,:
\begin{equation*}
\begin{aligned}
N^k:\Gr^W_{d_X+k}\M&\simto(\Gr^W_{d_X-k}\M)(-k),\\
N^k:\Gr^W_{d_X+1+k}V_{x,1}&\simto(\Gr^W_{d_X+1-k}V_{x,1})(-k),\\
N^k:\Gr^W_{d_X+k}V_{x,\ne1}&\simto(\Gr^W_{d_X-k}V_{x,\ne1})(-k).
\end{aligned}
\end{equation*}
These imply the $N$-{\it primitive decompositions}
\begin{equation*}
\begin{aligned}
\Gr^W_j\M&=\mopl_{k\ges0}\,N^k(P_N\Gr^W_{j+2k}\M)(k),\\
\Gr^W_jV_{x,1}&=\mopl_{k\ges0}\,N^k(P_N\Gr^W_{j+2k}V_{x,1})(k),\\
\Gr^W_jV_{x,\ne1}&=\mopl_{k\ges0}\,N^k(P_N\Gr^W_{j+2k}V_{x,\ne1})(k),
\end{aligned}
\end{equation*}
where the $N$-{\it primitive\one} parts $P_N\Gr^W_{\ssb}\M$, $P_N\Gr^W_{\ssb}V_{x,1}$, $P_N\Gr^W_{\ssb}V_{x,\ne1}$ are defined by
\begin{equation*}
\begin{aligned}
P_N\Gr^W_{d_X+k}\M&\defs{\rm Ker}\,N^{k+1}\sst\Gr^W_{d_X+k}\M,\\ 
P_N\Gr^W_{d_X+1+k}V_{x,1}&\defs{\rm Ker}\,N^{k+1}\sst\Gr^W_{d_X+1+k}V_{x,1},\\ P_N\Gr^W_{d_X+k}V_{x,\ne1}&\defs{\rm Ker}\,N^{k+1}\sst\Gr^W_{d_X+k}V_{x,\ne1},\end{aligned}
\end{equation*}
if $k\gess0$, and they vanish for $k\slt0$. It is quite well known that
\begin{equation*}
P_N\Gr^W_{d_X+k}\M\eq\begin{cases}\msum_x\,P_N\Gr^W_{d_X+k}V_{x,1}&\h{if}\,\,\,k\sgt0,\\ {\rm IC}_X\R_h&\h{if}\,\,\,k\eq0,\end{cases}
\end{equation*}
and it vanishes otherwise, see for instance \cite[Remark 2.3a]{FPS}, \cite[(2.2.5)]{KLS}. Here mixed Hodge modules supported on $x$ are identified with mixed Hodge structures and the direct image $(i_x)_*$ is omitted. We have also the isomorphism
\begin{equation*}
\Gr^W_{d_X}(\R_{h,X}[d_X])={\rm IC}_X\R_h,
\end{equation*}
(see \cite[(4.5.9)]{mhm}) and the vanishing cycle exact sequence of mixed Hodge modules
\begin{equation*}
0\tos\R_{h,X}[d_X]\tos\psi_g\R_{h,\Yt}[d_X]\tos\varphi_g\R_{h,\Yt}[d_X]\tos0,
\end{equation*}
as well as its restriction to the unipotent monodromy part $\psi_{g,1}$, $\varphi_{g,1}$. These imply the isomorphisms of mixed Hodge structures
\begin{equation*}
\begin{aligned}\mopl_{x\in{\rm Sing}\,X}H^{d_X-1}(L_{X,x},\R)&=C\bl(\R_{h,X}[d_X]\onto{\rm IC}_X\R_h\br)[-1]\\&={\rm Ker}\bl(\R_{h,X}[d_X]\onto{\rm IC}_X\R_h\br)\\&=W_{d_X-1}(\R_{h,X}[d_X])\\&={\rm Ker}\bl(N\,|\,\varphi_{g,1}\R_{h,\Yt}(1)[d_X]\br)\\&=\mopl_{x\in{\rm Sing}\,X}\,{\rm Ker}\bl(N\col V_{x,1}(1)\to V_{x,1}\br).\end{aligned}
\end{equation*}
Here the filtration on $\M$ defined by ${\rm Ker}\,N^j$ is used for the proof of the fourth isomorphism. (Note that ${\rm Ker}\,N$ commutes with $\Gr^W_{\ssb}$ by the strictness of morphisms of mixed Hodge modules.) We thus get the isomorphism \eqref{3}.
\sk
Since the two intersection complex $L$-classes coincides by Proposition\,\,\ref{P3.4}, the remaining assertions follow using \cite[Corollary 1.1]{Sc0}. Indeed, the above short exact sequence implies the second equality of \eqref{5} by using $\si(H(m))\eq(-1)^m\si(H)$ for $m\ins\Z$. (A similar argument related to the Milnor class is employed for instance in \cite{CMSS1}, \cite{MSS0}, \cite{MSS2}, using the virtual Hirzebruch class $T_{y*}^{\rm vir}(X)$.) The first equalities of \eqref{4} and \eqref{5} follow from the restriction of the above exact sequence to the unipotent monodromy part together with the isomorphism $\Gr^W_{d_X}(\R_{h,X}[d_X])\eq{\rm IC}_X\R_h$ and the hard Lefschetz property of the action of $N$ explained above.
\sk
As for the second equality of \eqref{4}, some more explanations is required.
The key point is the commutativity with the Gysin morphism and others under closed immersions, see \cite[Corollary 1.1]{Sc0} (and also \cite[Proposition 3.3]{MSS0}), that is,
\begin{equation*}
T_{1*}\bl(\psi_{g-c}\R_{h,\Yt}[d_X]\br)\eq i_c^!T_{1*}\bl(\R_{h,\Yt}[d_{\Yt}]\br)\eq i_c^!L_*(\Yt).
\end{equation*}
where $i_c\col X_c\eq g^{-1}(c)\into\Yt$ is the inclusion for $c\ins\ga$ (including $c\eq0$) and $g$ is identified with a function on the complement of $g^{-1}(\infty)$. Here we need also the last commutative diagram in \cite[Section 2.3]{MSS0} showing that
\begin{equation*}
i_c^!L_*(\Yt)=i_c^!\bl(L^*(\Yt)\cap[\Yt]\br)=i_c^*L^*(\Yt)\cap[X_c],
\end{equation*}
where $L_*(\Yt)\eq T_{1*}(\Yt)$, $L^*(\Yt)\eq T_1^*(\Yt)$, and $[X]$ denotes the fundamental homology class $[X]_{\rm hom}$. We have $i_c^*L^*(\Yt)\eq L^*(X_c)$ for $c\ins\ga\stm\{0\}$ and
\begin{equation*}
i_0^*L^*(\Yt)\eq L^*_{\rm vir}(X)\defs T^*_1(T_{\rm vir}X),
\end{equation*}
(since $N_{X_c/\Yt}$ is trivial for any $c\ins\ga$), hence
\begin{equation*}
i_0^!L_*(\Yt)=L^*_{\rm vir}(X)\cap[X]=T_{1*}^{\rm vir}(X)=:L_*^{\rm vir}(X),
\end{equation*}
where the latter is called the virtual $L$-class of $X$, see for instance \cite[Section 1.3]{MSS0} for $T^*_y(T_{\rm vir}X)$, $T_{y*}^{\rm vir}(X)$. We see that
\begin{equation*}
r^*L^*_{\rm vir}(X)=L^*(Z),
\end{equation*}
considering the restriction of $L^*(\Yt)$ to $g^{-1}(\ga)$, which is sent to $L^*_{\rm vir}(X)$ and $L^*(Z)$ under the restriction to $X_0$ and $Z$ respectively (and using a formula for $r^*$ noted at the beginning of this section). By the projection formula for $r$, we then get
\begin{equation*}
r_*L_*(Z)=r_*\bl(r^*L^*_{\rm vir}(X)\cap[Z]\br)=L^*_{\rm vir}(X)\cap r_*[Z]=L_*^{\rm vir}(X),
\end{equation*}
since $r_*[Z]\eq[X]$. So the assertion follows. This finishes the proof of Theorem\,\,\ref{T3}.

\begin{rem}\label{R8.1}
A good retraction $r\col Z\to X$ over $\ga$ can be obtained by using the integral curves of a $C^{\infty}$ real vector field $\xi$ on $g^{-1}(\ga)$ constructed by using a partition of unity so that $\langle\xi,g^*\ddd s\rangle\eq g^*s$ as in a proof of Thom's second isotopy lemma in the submersion (that is, trivial Whitney stratification) case. Here $s$ is a real coordinate of $\ga$ identifying it with $[0,1]$. Moreover $s$ is locally the restriction of a complex analytic local coordinate $z$ of $\C$ around 0, and $\ga$ is locally contained in $\{{\rm Re}\,z\gess0,\,{\rm Im}\,z\eq0\}$. We use also an embedded resolution $\pi$ near the singular fiber $X$ together with the action of $[0,1]$ constructed as in \cite{Cl}, and take the composition with $\pi$. By construction the retraction $r$ can be extended to a continuous mapping
\begin{equation*}
\rt:Z{\times}[0,1]\to g^{-1}(\ga),
\end{equation*}
such that $\rt|_{Z\times\{0\}}\eq r$, $\rt|_{Z\times\{1\}}\eq {\rm id}$, and $\rt|_{Z\times\{c\}}$ is a homeomorphism $Z\simto g^{-1}(s^{-1}(c))$ for $c\ins(0,1)$, where $\ga$ is identified with $[0,1]$ by $s$. We then see that $g^{-1}(\ga)$ is identified with the quotient topological space $(Z{\times}[0,1])/{\sim}$ (endowed with the quotient topology) by the relation: $(z,c)\sim(z',c')$ if and only if $r(z)\eq r(z')$ and $c\eq c'\eq0$. (Indeed, $Z{\times}[0,1]$ is compact and the quotient map is proper). So the composition noted at the beginning of this section coincides with $r^*$ considering the continuous map $g^{-1}(\ga)\cong(Z{\times}[0,1])/{\sim}\tos X{\times}[0,1]$ and the induced morphism between the higher direct images of the constant sheaves on these spaces under the projection to $[0,1]$.
\sk
The restriction of $r$ over the complement of ${\rm Sing}\,X$ is a homeomorphism, since $X$ is reduced (otherwise $X$ has non-isolated singularities). We then get the equality $r_*[Z]\eq[X]$, since we have an isomorphism $H^{\rm BM}_{2n}(X,\Z)\simto H^{\rm BM}_{2n}(X\stm{\rm Sing}\,X,\Z)$ induced by the open restriction. Note that the Bore-Moore homology fundamental class can be defined by using the canonical morphism $\Z_X[2d_X]\tos\DD\Z_X$, where $\DD\Z_X$ is the (topological) dualizing complex (and the Tate twist is omitted choosing $\sqrt{-1}$).
\end{rem}

\section{Proof of Proposition~\ref{P2}} \label{S9}
Set $X\defs\{f\eq0\}\sst\PP^3$, where
\begin{equation*}
f\defs x^{2a}z^b\pl y^{2a'}w^{b'},
\end{equation*}
with $a,a',b,b'$ mutually prime positive odd integers satisfying $2a\pl b\eq2a'\pl b'$; for instance, $(a,a',b,b')\eq(3,5,11,7)$. By the Thom-Sebastiani type theorem (see for instance \cite{MSS1}) we see that
\begin{equation}
C(\R_{h,X}[2]\tos{\rm IC}_X\R_h)[-1]=j_!L_U[1],
\label{9.1} \end{equation}
where $L_U$ is a local system of rank 1 on $U\defs\{x\eq y\eq0\}\stm\{zw\eq0\}\cong\C^*$ with $j\col U\into X$ the inclusion. Indeed, the vanishing cohomology of $g\eq u^av^b$ with $a,b$ mutually prime is the same as that of $g\eq uv$, in particular, the Milnor monodromy is unipotent. The vanishing cohomology of $g\eq u^a$ has no unipotent monodromy part. We verify that the monodromy of the local system $L_U$ is given by multiplication by $-1$ (exchanging the two local branches), so all the cohomology groups of $j_!L_U$ vanish, and the intersection cohomology coincides with the constant coefficient cohomology. (Note that $x^{2a}\mi\la\one y^{2b}\eq(x^a\mi\la^{1/2}y^b)(x^a\pl\la^{1/2}y^b)$.)
\sk
We then see that the intersection complex is the direct image of the shifted constant sheaf on the normalization of $X$, hence $L_U$ underlies a variation of Hodge structure of type $(0,0)$, and we can apply the argument in Section\,\,\ref{S4} on the graded pieces of the topological filtration $G$. We also verify that the Hodge module whose underlying $\R$-complex is $j_!L_U[1]$ belongs to the image of Hdg in \eqref{1} considering the difference between the classes of $\overline{U}\cong\PP^1$ and a double covering of $\overline{U}$ ramified over $\overline{U}\stm U$. This finishes the proof of Proposition\,\,\ref{P2}.

\begin{rem} \label{R9.1}
The above example can be extended to the case $f\eq\msum_{i=1}^m\,x_i^{a_i}x_{i+m}^{b_i}$ ($m\gess 2$). Here, setting $a'_i\eq a_i/2$ for $i\eq1,2$, $a'_i\eq a_i$ for $i\sgt2$, and $b'\eq b_i$ for any $i$, the $2m$ numbers $a'_i,b'_i$ ($i\ins[1,m]$) are mutually prime odd integers, and the $a_i\pl b_i$ are independent of $i\ins[1,m]$.
\end{rem}

\section{Images of Hdg and sd} \label{S10}
It seems very difficult to prove or disprove that the classes of intersection complexes belong to the images of the morphisms $\rm Hdg$ and $sd$ in \eqref{1}. We first consider the case of the morphism $sd$, since this seems simpler. Let $\pi\col X'\tos X$ be a birational morphism with $X'$ smooth projective. The direct factors $^p\!R^0\pi_*(\R_{X'}[d_{X'}])_Z$ of $^p\!R^0\pi_*(\R_{X'}[d_{X'}])$ with strict supports $Z\sst X$ (see \cite[(5.1.3.5)]{mhp}) are intersection complexes with coefficients in geometric variations of Hodge structure, which are {\it direct factors\one} of $^p\!R^0\!\one g_*(\R_V[d_V])_Z$ with $V$ smooth projective and $g\col V\tos Z$ surjective. Here one can use the {\it weak Lefschetz\one}-type theorem in order to get $^p\!R^0\!\one g_*(\R_V[d_V])$, but this may imply a necessity of taking a direct factor, since we have only the {\it injectivity\one} at the relevant degree in that theorem. Here we study the image of $sd$ using for instance \cite[Proposition 1.1a]{FPS}, see also \cite{CS}, \cite{Yo}.
\sk
Consider the case where $X$ is a compactification of the (disjoint) union of the projective cones $C(V_s)$ of $n$-dimensional smooth projective varieties $V_s$ parametrized by a smooth variety $S$ and $X'$ is a smooth compactification of the blow-up of $\bigsqcup_{s\in S}\,C(V_s)$ along the union of vertices (which is identified with $S$), where $Z$ is a compactification of $S$. It is well known that the restriction of the direct factor $^p\!R^j\pi_*(\R_{X'}[d_X])_Z$ to $S$ is given by the {\it non-primitive\one} part of $H^{n+1+j}(V_s)$ (up to a shift of complex), since the primitive part gives the stalks of the intersection complex of the union of the cones using the Thom-Gysin sequence, see for instance \cite[1.3]{RSW}. Here we have a shift of the cohomological degree by 1, and this can be adjusted by taking {\it general hyperplane sections\one} $V'_s$ of the $V_s$ and applying the hard and weak Lefschetz theorems so that we get the inclusion
\begin{equation*}
H^{n+1}(V_s)\cong H^{n-1}(V_s)\into H^{n-1}(V'_s),
\end{equation*}
which splits using a polarization. However, this injection is never bijective in general.
\sk
It seems then interesting whether we can get a family $V'_s$ {\it without using the hyperplane sections.} Assume the $V_s$ are abelian surfaces and the $V'_s$ are curves of {\it genus\one}~2. In this case every curve is hyperelliptic and the moduli spaces of curves and their Jacobians have the same dimension 3, but the theta divisor, which coincides with the image of the curve and is non-singular by Riemann (see for instance \cite[VI,\,Section 1]{ACGH}), is not {\it very\one} ample, since we would get a degenerating family of curves otherwise. Here one can replace $S$ with a non-empty open subvariety as long as one proceeds by induction on $\dim V_s\pl\dim S$ (since the support of the unipotent monodromy part of the vanishing cycle complexes in the normal crossing case has codimension 2, see also \cite[Corollary 4.2.4]{mhp}).
\sk
It seems completely unclear whether the answer is positive or not because of the difference between {\it fine\one} and {\it coarse\one} moduli spaces. The Torelli theorem does not seem enough to solve a problem related to {\it fine\one} moduli spaces because of automorphisms of curves, and the involution of hyperelliptic curves may cause a problem. Using spreading-out, one can construct a family $V'_s$ over a {\it finite covering\one} of $S$ after shrinking $S$, but this is not sufficient to solve the problem, since we have to take direct factors of direct images as before.

\begin{rem} \label{R10.1}
The above problem is closely related to the one about the image of Hdg in \eqref{1}. Assume for instance there is a family of genus 2 curves $C_s$ on a smooth variety $S$ such that the $\R$-local system $L'$ of the first cohomology groups is isomorphic to that of a family of abelian surfaces $V_s$ over $S$, which is denoted by $L$ (shrinking $S$ if necessary). Assume the family $\{V_s\}_{s\in S}$ has {\it entirely non-trivial variation,} that is, $L_{\C}\cap F^1\Lc\eq0$, where $(\Lc,F)$ is the underlying filtered $\OO_S$-module of the corresponding variation of Hodge structure of fiber-wise weight 1 with $\Lc\eq\OO_S{\otimes}_{\R}L$, $L_{\C}\eq\C{\otimes}_{\R}L$. Using a global invariant cycle theorem (which follows from the adjunction morphism $a_S^*(a_S)_*\tos id$ with $a_S\col S\tos pt$ the structure morphism), we see that the isomorphism is {\it compatible\one} with the Hodge filtration of the variations of Hodge structure. (Indeed, the isomorphism gives a global section of the local system $L''\defs{\Hc}\!om(L',L)$, and its $(1,-1)$-component is also a global section of $L''$ by that theorem, but its image must vanish by the assumption on entirely non-trivial variation.)
\sk
The problem is thus closely related to the question whether the classes of intersection complexes belong to the image of the morphism Hdg in \eqref{1}. Here one can construct an isomorphism of $\Q$-local systems, but not $\Z$-local systems. So the theory related to non-principally polarized abelian varieties cannot be used. (See \cite[Proposition 1]{FPS} for the independence of polarizations over $\R$.)
\end{rem}

\begin{rem} \label{R10.2}
The above problems seem quite non-trivial even in the case where the $V_s$ are self-products of elliptic curves $E_s$. Here it seems interesting whether the $V'_s$ may be a hyperelliptic curves of genus 2 having an involution such that the quotient of $V'_s$ by this involution is $E_s$. Let $E_s\tos C\defs\PP^1$ be a finite covering of degree 2 ramified over $0,1,\la,\infty\ins C$. Let $C'_u\defs\PP^1\to C$ be a finite covering of degree 2 ramified over $u,\infty\ins C$. Let $E'_{s,u}$ be the normalization of $E_s{\times}_CC'_u$ with $\rho\col E'_{s,u}\tos E_s$ the canonical morphism of degree 2. Here the unique singular point of $E_s{\times}_CC'_u$ is over $\infty\ins C$, and has normal crossings, since it is locally expressed as $\{z^2\eq w^2\}$ by the definition of the fiber product. We then see that $\rho$ is an \'etale morphism except over the two points $e',e''\ins E_s$ consisting of the inverse image of $e\ins C$.
\sk
We have a canonical injection $H^1(E_s,\Q)\into H^1(E'_{s,u},\Q)$, and its cokernel is a Hodge structure of weight 1 and level 1, which depends on $u$, and would be isomorphic to $H^1(E_s,\Q)$ for some $u_s$ depending on $s$. (Here we use $j$-invariant, see \cite{Si}.) But the problem is whether we can get a {\it univalent\one} branch when $s$ varies globally. Note that $\la$ depends on $s$, and is not necessarily {\it univalent.} Assume the elliptic curve $E_s$ is defined for instance by the equation
\begin{equation*}
y^2\eq x^3\pl sx\pl s\q(s\ins S\defs\C\stm\{0,-\tfrac{27}{4}\}).
\end{equation*}
Its $j$-invariant is given by $s/(4s\pl27)$ up to a non-zero constant multiple, see \cite[p.\,50]{Si}. This may suggest that we have a {\it unique\one} choice of $u_s$ for general $s\ins S$, or more precisely, there is a univalent section on an open subvariety of $S$. However, it seems quite difficult to prove its existence or non-existence.
\sk
It may be interesting to consider also the case where the $V_s$ are products of elliptic curves and $S$ is the self-product of $\C\stm\{0,-\tfrac{27}{4}\}$ (or the quotient of the complement of its diagonal by the action of ${\mathfrak S}_2$). Here it does not seem easy to control finite morphisms from curves of genus 2 to elliptic curves, since its degree is not necessarily 2. (If it is 2, we can get an involution of the curve.) For instance, there is an example where the degree is 3 and it is ramified over only one point. Indeed, there are $\si,\si'\ins{\mathfrak S}_3$ with order 2 (which are the monodromies along two loops on the elliptic curve) such that $\si'{}^{-1}\si^{-1}\si'\si$ (which is the inverse of the local monodromy at the ramification point) has order 3; for instance $\si\eq(1,2)$, $\si'\eq(2,3)$.
\end{rem}

\begin{rem} \label{R10.3}
It is quite well known that any polarizable $\R$-Hodge structure is determined only by Hodge numbers up to non-canonical isomorphisms (more precisely, it is a direct sum of $\R$-Hodge structures of rank 1 or 2), and is ``geometric" so that its class belongs to the image of Hdg in \eqref{1} with $X$ a point. (This can be verified by increasing induction on the level of Hodge structure using self-products of an elliptic curve for instance.)
\sk
By an argument similar to \cite[(5.3.3)]{MSS0}, we can employ the above assertion to prove that the class of the intersection complex $\R$-Hodge module $\bl[{\rm IC}_X\R_h\br]$ of a compact variety $X$ belongs to the image of Hdg in \eqref{1} if the restriction of the cohomology sheaves $\Hc^j{\rm IC}_X\R$ to $S$ are ``constant" for any $j\ins\Z$ and any stratum $S$ of an algebraic stratification of $X$. (It is unnecessary to assume that it is a Whitney stratification, since we can take a refinement.) To show that the variations of Hodge structure on strata are ``constant", one can apply a functorial morphism $(a_S)_*\tos i_s^*$ induced from $id\tos(i_s)_*i_s^*$, where $i_s\col\{s\}\tos S$ and $a_S\col S\tos pt$ are natural morphisms. (This gives a generalization of the global invariant cycle theorem.)
\end{rem}

\bs

\end{document}